\documentclass[reqno]{amsart}

\usepackage{xypic}
\input xy
\xyoption{all}
\usepackage{epsfig}
\usepackage{color}
\usepackage{amsthm}
\usepackage{amssymb}
\usepackage{amsmath}
\usepackage{amscd}
\usepackage{amsopn}

\newcommand{\labell}[1] {\label{#1}}

\numberwithin{equation}{section}
\newtheorem {Theorem}{Theorem}
\numberwithin{Theorem}{section}

\newtheorem {Lemma}[Theorem]    {Lemma}
\newtheorem {Question}[Theorem]    {Question}
\newtheorem {Proposition}[Theorem]{Proposition}
\newtheorem {Corollary}[Theorem]{Corollary}
\theoremstyle{definition}
\newtheorem{Definition}[Theorem]{Definition}
\theoremstyle{remark}
\newtheorem{Remark}[Theorem]{Remark}
\newtheorem{Example}[Theorem]{Example}

\expandafter\chardef\csname pre amssym.def at\endcsname=\the\catcode`\@
\catcode`\@=11
\def\undefine#1{\let#1\undefined}
\def\newsymbol#1#2#3#4#5{\let\next@\relax
 \ifnum#2=\@ne\let\next@\msafam@\else
 \ifnum#2=\tw@\let\next@\msbfam@\fi\fi
 \mathchardef#1="#3\next@#4#5}
\def\mathhexbox@#1#2#3{\relax
 \ifmmode\mathpalette{}{\m@th\mathchar"#1#2#3}%
 \else\leavevmode\hbox{$\m@th\mathchar"#1#2#3$}\fi}
\def\hexnumber@#1{\ifcase#1 0\or 1\or 2\or 3\or 4\or 5\or 6\or 7\or 8\or
 9\or A\or B\or C\or D\or E\or F\fi}

\font\teneufm=eufm10
\font\seveneufm=eufm7
\font\fiveeufm=eufm5
\newfam\eufmfam
\textfont\eufmfam=\teneufm
\scriptfont\eufmfam=\seveneufm
\scriptscriptfont\eufmfam=\fiveeufm

\catcode`\@=\csname pre amssym.def at\endcsname

\def    \eps    {\epsilon}

\newcommand{\FF}{{\mathcal F}}

\newcommand{\CC}{{\mathcal C}}

\newcommand{\CM}{{\mathcal M}}
\newcommand{\CN}{{\mathcal N}}
\newcommand{\CU}{{\mathcal U}}

\newcommand{\CS}{{\mathcal S}}

\newcommand{\supp}{\operatorname{supp}}

\newcommand{\Ham}{{\mathit Ham}}

\newcommand{\id}{{\mathit id}}
\newcommand{\const}{{\mathit const}}

\newcommand{\CB}{{\mathcal B}}

\newcommand{\PP}{{\mathcal P}}

\newcommand{\CW}{{\mathcal W}}

\def    \R      {{\mathbb R}}

\def    \Z      {{\mathbb Z}}

\def    \T      {{\mathbb T}}

\def    \12    {{\frac{1}{2}}}
\def    \bG    {\bar{\Gamma}}
\def    \p      {\partial}

\def    \length  {\operatorname{length}}

\def    \ev  {\operatorname{ev}}

\def    \HF     {\operatorname{HF}}
\def    \HM     {\operatorname{HM}}

\def    \CF     {\operatorname{CF}}

\def    \Ham     {\operatorname{Ham}}
\def    \Fix     {\operatorname{Fix}}

\def    \MUCZ  {\operatorname{\mu_{\scriptscriptstyle{CZ}}}}

\def    \ssminus        {\smallsetminus}

\def    \chom {\operatorname{c_{hom}}}

\begin{document}


\setlength{\smallskipamount}{6pt}
\setlength{\medskipamount}{10pt}
\setlength{\bigskipamount}{16pt}





\title[Coisotropic Intersections]{Coisotropic Intersections}

\author[Viktor Ginzburg]{Viktor L. Ginzburg}

\address{Department of Mathematics, UC Santa Cruz,
Santa Cruz, CA 95064, USA}
\email{ginzburg@math.ucsc.edu}

\date{\today}
\thanks{The work is partially supported by the NSF and by the faculty
research funds of the University of California, Santa Cruz.}

\bigskip

\begin{abstract}
In this paper we make the first steps towards developing a theory of
intersections of coisotropic submanifolds, similar to that for
Lagrangian submanifolds.

For coisotropic submanifolds satisfying a certain stability
requirement we establish persistence of coisotropic intersections
under Hamiltonian diffeomorphisms, akin to the Lagrangian intersection
property. To be more specific, we prove that the
displacement energy of a stable coisotropic submanifold is positive,
provided that the ambient symplectic manifold meets some natural
conditions.  We also show that a displaceable, stable, coisotropic
submanifold has non-zero Liouville class. This result further
underlines the analogy between displacement properties of Lagrangian
and coisotropic submanifolds.

\end{abstract}

\maketitle

\section{Introduction}
\labell{sec:intro} In this paper we make the first steps towards
developing a theory of coisotropic intersections similar to that for
Lagrangian submanifolds. The main objective of the paper is to
establish persistence of coisotropic intersections under Hamiltonian
diffeomorphisms for a certain class of coisotropic submanifolds,
the so-called stable coisotropic submanifolds.  We also show that
displaceable, stable coisotropic submanifolds have  non-zero Liouville
class, generalizing the results of Bolle, \cite{Bo1,Bo2}, for
submanifolds of $\R^{2n}$.

\subsection{Coisotropic intersections}
The Lagrangian intersection property or persistence of Lagrangian
intersections is unquestionably one of the most fundamental results in
symplectic topology. This result asserts that a Lagrangian submanifold
necessarily intersects its image under a Hamiltonian diffeomorphism
that is in some sense close to the identity, e.g., has sufficiently
small energy.  Depending on the notion of closeness and on the
requirements imposed on the manifolds, various forms of the
Lagrangian intersection property have been proved in
\cite{Ch:Cambr,Ch,F:Morse,F:grad,F:c-l,Fl,Gr,LS:lagr,Oh1,Oh:cambridge,Oh:disj,We73,We74},
to mention just some of the pertinent references.

Recall that a submanifold $M$ of a symplectic manifold
$(W^{2n},\omega)$ is said to be coisotropic if for every $p\in M$
the symplectic orthogonal $(T_pM)^\omega$ to the tangent space
$T_pM$ is contained in $T_pM$.  For instance, Lagrangian
submanifolds are coisotropic, as are hypersurfaces in $W$.
Furthermore, $\dim M\geq n$, when $M$ is coisotropic. The examples
discussed below indicate that coisotropic submanifolds enjoy the
same kind of Hamiltonian rigidity as Lagrangian submanifolds and
lead to the following

\begin{Question}[Coisotropic intersections]
\label{quest:intersections} Can a coisotropic submanifold be
displaced by a Hamiltonian diffeomorphism arbitrarily close to the
identity in a suitable sense? For instance, interpreting closeness
in the sense of Hofer's metric (see, e.g.,
\cite{HZ,Pol:book}), we can ask if there is a lower bound, depending
only on the submanifold, on the energy of a Hamiltonian
diffeomorphism displacing the submanifold.
\end{Question}

This question can be further restricted in a number of ways: by
imposing additional assumptions on the codimension of $M$ or on the
dynamics of the characteristic foliation or via other types of
symplectic-topological requirements on $M$ and the ambient manifold,
or through requiring the Hamiltonian diffeomorphism to be close to
the identity in a particularly strong way. For instance, when $M$ is
Lagrangian, the question reduces to the Lagrangian intersection
property. In terms of codimension, the other extreme case is that of
codimension zero: $M=W$. In this case, the coisotropic intersection
property obviously holds. When $M$ is a hypersurface, the answer to
the question is affirmative due to non-degeneracy of Hofer's metric,
\cite{LMD95}, and the fact that a displaceable connected
hypersurface necessarily bounds. The extreme case in terms of
closeness of the diffeomorphism to the identity is that of
infinitesimal displacement:

\begin{Example}[Infinitesimal intersections]
\label{ex:infinitesimal} Let $M$ be a coisotropic submanifold and
let $H$ be a Hamiltonian near $M$. On the infinitesimal level, the
intersections of $M$ with its image under the Hamiltonian flow of
$H$ correspond to the points where the Hamiltonian vector field
$X_H$ of $H$ is tangent to $M$. These are precisely the leaf-wise
critical points of $H$ along the characteristic foliation $\FF$ on
$M$. This observation suggests that the ``amount'' of coisotropic
intersections is governed by the foliated Morse theory of $\FF$.
\end{Example}

When $M$ is Lagrangian, this observation readily implies persistence
of intersections for any Hamiltonian diffeomorphism generated by a
$C^2$-small Hamiltonian, by the Weinstein symplectic neighborhood
theorem, \cite{We71,We77}.

Here we answer Question \ref{quest:intersections} affirmatively only
for coisotropic manifolds satisfying a certain additional stability
requirement (see Definition \ref{def:B}), introduced by Bolle in
\cite{Bo2} and similar to stability of hypersurfaces (cf.\
\cite{EKP,HZ}).  Namely, for such a submanifold we establish a lower
bound on the Hofer norm of a displacing Hamiltonian diffeomorphism
(see Theorem \ref{thm:main}) and hence prove that the submanifold has
positive displacement energy.  The stability condition, discussed in
detail in Section \ref{sec:ct}, is quite restrictive; see, in
particular, Example \ref{exam:Bolle}. (We also impose some natural and
not-too-restrictive assumptions on the ambient manifold.)
Nevertheless, these results combined with the examples given in this
section appear to provide enough evidence to conjecture that the
coisotropic intersection property holds in general.

Next let us examine the case where the characteristic foliation
is a fibration. We do this by passing to the graph of the
foliation and interpreting it as a Lagrangian submanifold. This will also lead
us to a refinement of Question \ref{quest:intersections}.

\begin{Example}[Leaf-wise coisotropic intersections]
\label{ex:leafwise}
Let $M$ be a coisotropic submanifold of $(W,\omega)$ and let
$\Gamma\subset M\times M\subset W\times W$ be the graph of the
characteristic foliation $\FF$ on $M$. In other words, $\Gamma$ is
formed by pairs $(x,y)\in M\times M$ with $x$ and $y$ lying on the
same leaf of $\FF$. Then $\Gamma$ is a one-to-one immersed Lagrangian
submanifold of $W\times W$, where the latter is equipped with the
symplectic form $\omega\oplus -\omega$. (In general, $\Gamma$ is not a
true submanifold, e.g., $\Gamma$ can be dense in $M\times M$.)
Consider a Hamiltonian diffeomorphism $\varphi$ of $W$. Then
$\tilde{\varphi}=(\id,\varphi)$ is a Hamiltonian diffeomorphism of
$W\times W$ and the intersection points
$\tilde{\varphi}(\Gamma)\cap\Gamma$ are in one-to-one correspondence
with $x\in M$ such that $x$ and $\varphi(x)\in M$ lie on the same
leaf. Hence, persistence of Lagrangian intersections for $\Gamma$, if
it held, would imply the \emph{leaf-wise intersection property} for
$M$, i.e., the existence of a leaf $F$ of $\FF$ with $\varphi(F)\cap
F\neq \emptyset$. For instance, assume that $\FF$ is a fibration. Then
$\Gamma$ is a true smooth Lagrangian submanifold of $W\times W$ and we
conclude that in this case leaf-wise intersections do exist, provided
that $\varphi$ is not far from $\id$ in a suitable sense; see, e.g.,
\cite{Ch,F:Morse,F:grad,F:c-l,Fl,Gr,LS:lagr,Oh1,Oh:disj,We73}.
\end{Example}

Taking this example as a motivation let us call $x\in M$ a leaf-wise
intersection of $M$ and $\varphi(M)$ if $x\in F\cap\varphi(F)$ for some
leaf $F$ of the coisotropic foliation. Note that in the infinitesimal
setting of Example \ref{ex:infinitesimal}, leaf-wise intersections
correspond to the critical points of $H$ on $M$.

\begin{Question}[Leaf-wise coisotropic intersections]
\label{quest:lw-intersections} Do leaf-wise intersections exist
(perhaps, under some additional conditions on the coisotropic
submanifold) whenever the Hamiltonian diffeomorphism is
sufficiently close to the identity, e.g., in the sense of Hofer's
metric?
\end{Question}

When $M$ is Lagrangian, this question is of course equivalent to
Question \ref{quest:intersections}, for the characteristic foliation
in this case has only one leaf, the entire manifold $M$. In the other
extreme case $M=W$, where $W$ is closed, Question
\ref{quest:lw-intersections} is equivalent to the Arnold conjecture --
the existence of fixed points of Hamiltonian diffeomorphisms -- and
hence the answer to the question is affirmative in this case. (See,
e.g., \cite{F:Morse,F:grad,F:c-l,Fl,FO,Gr,HZ,LT:ac,mdsa-book,Sa} and references
therein.)

Question \ref{quest:lw-intersections}, arising also from some problems
in Hamiltonian dynamics, was originally raised by Moser in 1978 in
\cite{Mo}.  In \cite{Ba,Mo}, persistence of leaf-wise intersections
was proved for closed coisotropic submanifolds and Hamiltonian
diffeomorphisms which are $C^1$-close to the identity.  Moser's theorem
was later extended by Hofer,
\cite{Ho90}, to hypersurfaces of restricted contact type in $\R^{2n}$
and Hamiltonian diffeomorphisms with energy smaller than a certain
symplectic capacity of the region bounded by the
hypersurface. (See also \cite{EH,Li}.)
Recently, Dragnev, \cite{Dr}, generalized this result to
arbitrary closed coisotropic submanifolds in $\R^{2n}$ that have
contact type, but need not be of restricted contact type
(see Definition \ref{def:B}) . In this case, the energy of the
diffeomorphism is required to be smaller than the homological
capacity of the submanifold itself. Note that this capacity is
positive for any coisotropic submanifold which is displaceable and stable
as is shown below in Remark \ref{rmk:Dr}.

Here we do not consider the problem of extending Dragnev's theorem
to other ambient manifolds; this question will be addressed
elsewhere. Instead, we prove a simple preliminary result in this
direction and generalize the results of \cite{EH,Ho90} to
subcritical Stein manifolds.  Namely, we show that leaf-wise
intersections of hypersurfaces in such manifolds persist for
Hamiltonian diffeomorphisms with energy smaller than the
homological capacity of the region bounded by the hypersurface. (See
Theorem \ref{thm:leafwise}.) As
stated, this result does not hold when the contact type condition is
dropped; see Example \ref{ex:leafwise-fail}.

Finally note that Moser's theorem and the results on
Lagrangian intersections and the Arnold conjecture discussed above
suggest the existence, under suitable additional assumptions, of more
than one leaf-wise intersection.

\begin{Remark}[Totally non-coisotropic displacement, \cite{gu:new}]
A theorem of Polterovich and of Laudenbach and Sikorav,
\cite{LS,Pol}, asserts that Hamiltonian persistence of intersections
is an exclusive feature of Lagrangian submanifolds among
submanifolds of middle dimension.  Namely, a middle-dimensional
submanifold $N$ admits an infinitesimal Hamiltonian displacement if
and only if $N$ is not Lagrangian and its normal bundle has a
non-vanishing section.  This clear-cut dichotomy does not carry over
to lower codimensions, but the general picture is somewhat similar.
To be more specific, it has been shown by G\"urel, \cite{gu:new},
that a totally non-coisotropic submanifold admits an infinitesimal
Hamiltonian displacement, provided that its normal bundle has a
non-vanishing section. Now, in contrast with the Lagrangian
case, it is not sufficient to assume (under the same normal bundle
condition) that $N$ is simply not coisotropic, for $N$ may contain a
Lagrangian submanifold.
\end{Remark}

\subsection{Coisotropic Liouville class}
The Lagrangian intersection property is intimately connected to the
fact that a Lagrangian submanifold that is displaceable must have 
non-zero Liouville class. Moreover, under suitable hypotheses, the
displacement energy can be bounded from below via the size
of the Liouville class. This connection can also be extended to
coisotropic manifolds.

Thus assume, for the sake of simplicity, that the ambient symplectic
manifold $(W^{2n},\omega)$ is symplectically aspherical and exact:
$\omega=d\lambda$. Then the restriction $\lambda|_{\FF}$ of $\lambda$
to the characteristic foliation $\FF$ of $M$
is leaf-wise closed, and the cohomology class
$[\lambda|_{\FF}]\in H^1_{dR}(\FF)$ in the foliated de Rham cohomology
is defined.  (Recall that $H^*_{dR}(\FF)$ is the cohomology
of the complex of smooth differential forms along the leaves of $\FF$;
see, e.g., \cite{MS} and references therein for a
discussion of foliated de Rham cohomology.)  Note that $[\lambda|_{\FF}]$
depends on the choice of $\lambda$.  When $M$ is Lagrangian, this is
the ordinary Liouville class of $M$. In general, we will refer to it
as the \emph{coisotropic Liouville class}. By analogy with Lagrangian
manifolds, we ask

\begin{Question}[Coisotropic Liouville class]
\label{quest:Liouville}
Is $[\lambda|_{\FF}]\neq 0$, provided that $M$ is displaceable and
closed?
\end{Question}

When $M$ is Lagrangian, this answer is affirmative by
\cite{Ch:Cambr,Ch,Gr,Oh:disj,Po93}. If $M$ is a hypersurface, the
answer is also affirmative. (Indeed, a displaceable hypersurface
bounds a region and then $\int_M\lambda\wedge\omega^{n-1}$ is equal
to the symplectic volume of the region. On the other hand, this
integral would be zero if $\lambda$ were leaf-wise exact.)
Furthermore, again by analogy with the Lagrangian case, one can ask
whether the displacement energy of $M$ can be bounded from below via
the ``size'' of $[\lambda|_{\FF}]$. A sufficiently good lower
bound would imply an affirmative answer to Question
\ref{quest:intersections} by a version of the figure-eight trick.

\begin{Example}
\label{ex:Liouville-graph}
In the setting of Example \ref{ex:leafwise}, assume that
$\omega=d\lambda$.  Then the manifold $W\times W$ is exact and the
Liouville class $[\lambda_\Gamma]$ of the graph $\Gamma$ is defined,
when $\Gamma$ is interpreted as an immersed manifold. Denote by
$\pi_1$ and $\pi_2$ the natural projections of $\Gamma$ to the first
and the second factor in $W\times W$. Then
$\lambda_\Gamma=\pi_1^*\lambda-\pi_2^*\lambda$.  Furthermore, it is
not hard to show that $[\lambda|_{\FF}]\neq 0$ if
$[\lambda_\Gamma]\neq 0$. Assume now that $\FF$ is a fibration and $M$
is displaceable.  Then $\Gamma$ is a genuine, displaceable Lagrangian
submanifold of $W\times W$. Under natural additional
assumptions on $W$, we infer that $[\lambda_\Gamma]\neq 0$ (see, e.g.,
\cite{Ch:Cambr,Ch,Gr,Oh:disj,Po93}), and hence $[\lambda|_{\FF}]\neq
0$. Note also that, as a consequence, the first cohomology of the
fiber of $\FF$ is non-trivial.
\end{Example}

We prove that, for a stable coisotropic submanifold, there exists a
loop $\eta$ which is tangent to $\FF$, contractible in $W$ and bounds
non-zero symplectic area (Theorem~\ref{thm:main}(ii)). This theorem
requires natural minor assumptions on the behavior of $\omega$ at
infinity in $W$, but holds even when $W$ is not exact.  As a
consequence of the theorem, we obtain an affirmative answer to
Question \ref{quest:Liouville} for stable submanifolds. However, this
approach does not lead to an answer to this question in
general, for such a loop $\eta$ need not exist when the stability
condition is dropped even if the Liouville class is non-zero. (This
follows from counterexamples to the Hamiltonian Seifert
conjecture; see, e.g., \cite{Ci,Gi95,gi:bayarea,gg:ex2,ke:example} and
references therein.)

Furthermore, if $M$ has restricted contact type, the displacement
energy of $M$ is greater than or equal to the symplectic area bounded
by $\eta$. (See Definition \ref{def:B} and Theorem
\ref{thm:main}(iii).)

For $W=\R^{2n}$, these results were proved by Bolle, \cite{Bo1,Bo2},
using the finite--dimensional reduction methods of \cite{CZ83}. Our
proof draws heavily on Bolle's ideas and is in fact just a translation
of his argument to the Floer theoretic setting.\footnote{The author is
grateful to Claude Viterbo for calling his attention to Bolle's papers,
\cite{Bo1,Bo2}, which played a crucial role in this work.}

\subsection{Dense existence of non-contractible loops}
\label{sec:dense} The stability condition and the proof of Theorem
\ref{thm:main} suggest that the loop $\eta$, tangent to a leaf of the
characteristic foliation but not contractible in it, can be viewed as
a generalization of a closed characteristic on a hypersurface; cf.\
\cite{Bo1,Bo2}.  Then assertions (ii) and (iii) of Theorem
\ref{thm:main} are interpreted as generalizations of the existence of
closed characteristics on stable hypersurfaces in $\R^{2n}$,
established in \cite{HZ:cap, HZ, St}. (See also \cite{gi:alan} for
further references.)

Continuing this analogy, consider a map $\vec{K}\colon W\to\R^k$
whose components are proper, Poisson--commuting Hamiltonians. Then
$M_a=\vec{K}^{-1}(a)$ is a coisotropic submanifold in $W$ whenever
$a\in\R^k$ is a regular value of $\vec{K}$. Assume, for the sake of
simplicity, that all levels $M_a$ are displaceable.  Then we prove
that for a subset $A\subset\R^k$ dense in the set of regular values of
$\vec{K}$, a level $M_a$, with $a\in A$, carries a loop tangent to the
characteristic foliation, contractible in $W$ but not in its leaf,
and bounding a positive symplectic area. This
result (Theorem \ref{thm:comm}) can be regarded as a generalization
of the dense existence theorem, \cite{fh,HZ}; see also
\cite{gi:alan}.  Here, as above, we need to impose some natural
additional conditions on $W$, but the coisotropic submanifolds $M_a$
need not be stable. This theorem leads to the question whether a
version of the almost existence theorem, \cite{HZ:cap, HZ, St}, for
such loops holds for commuting Hamiltonians.

\subsection{Organization of the paper}
The main results of the paper and the necessary definitions are stated
and further discussed in detail in Section \ref{sec:statements}.

The goal of Section \ref{sec:non-deg} is purely technical: here we set
conventions and recall relevant results concerning filtered Floer
homology, homotopy maps, and action selectors.

In Sections \ref{sec:nondeg} and \ref{sec:deg} we establish an
auxiliary result on which the proofs of our main theorems hinge. This
result, which may be of independent interest
(cf.\ \cite{En,KL,LM,MDS,Oh:chain}), asserts the existence of a Floer
connecting trajectory descending from a one-periodic orbit to the
maximum of a Hamiltonian and having energy bounded from above by the
displacement energy of the support of the Hamiltonian, provided that
the maximum is large enough and the Hamiltonian is ``slow'' near its
maximum. For non-degenerate Hamiltonians this result (Proposition
\ref{prop:orbits-nondeg}) is proved in Section \ref{sec:nondeg}. In
Section \ref{sec:deg} we deal with the degenerate case (Proposition
\ref{prop:orbits-deg}) and also discuss the space of finite energy
Floer trajectories.

Section \ref{sec:prfs} is devoted to the proofs of the main results of
the paper.  Here we establish Theorem \ref{thm:main} giving
affirmative answers to Questions \ref{quest:intersections} and
\ref{quest:Liouville} for stable, coisotropic submanifolds. In this
section we also prove a version of the dense existence theorem for
commuting Hamiltonians (Theorem \ref{thm:comm}).

Finally, in Section \ref{sec:leafwise} we prove persistence of
leaf-wise intersections for hypersurfaces of restricted contact type
in subcritical Stein manifolds (Theorem \ref{thm:leafwise}).

\subsection*{Acknowledgments.} The author is deeply grateful to
Yasha Eliashberg, Ba\c sak G\"urel, Ely Kerman, Felix Schlenk, and
Claude Viterbo for their numerous valuable remarks and suggestions.

\section{Displacement of stable coisotropic submanifolds}
\label{sec:statements}
\subsection{Contact type and stable coisotropic submanifolds}
\label{sec:ct}
Let $(W^{2n}, \omega)$ be a symplectic manifold and let $M\subset W$ be a
closed coisotropic submanifold of codimension $k$. Set
$\omega_0=\omega|_M$.  Then, as
is well known, the distribution $\ker \omega_0$ has dimension $k$ and
is integrable.  Denote by $\FF$ the characteristic foliation on $M$,
i.e., the $k$-dimensional foliation whose leaves are tangent to the
distribution $\ker \omega_0$.

\begin{Definition}
\label{def:B}
The coisotropic submanifold $M$ is said to be \emph{stable} if there exist
one-forms $\alpha_1,\ldots,\alpha_k$ on $M$ such that $\ker
d\alpha_i\supset \ker\omega_0$ for all $i=1,\ldots,k$ and
\begin{equation}
\labell{eq:ct}
\alpha_1\wedge\cdots\wedge\alpha_k\wedge \omega_0^{n-k}\neq 0
\end{equation}
anywhere on $M$. We say that $M$ has
\emph{contact type} if the forms $\alpha_i$ can be taken to be
primitives of $\omega_0$. Furthermore, $M$ has \emph{restricted}
contact type if the forms $\alpha_i$ extend to global primitives of
$\omega$ on $W$.
\end{Definition}

Stable and contact type coisotropic submanifolds were introduced by
Bolle in \cite{Bo1,Bo2}.  The nature of the requirements of Definition
\ref{def:B}, which are very restrictive, is illustrated by the
following examples.

\begin{Example}
\label{exam:Bolle}
~
\begin{enumerate}
\item[(i)] A contact type coisotropic submanifold is automatically
stable. Furthermore, a coisotropic submanifold which is $C^1$-close to
a coisotropic submanifold of contact type also has contact
type. (Apparently the latter is not true for stable coisotropic
submanifolds.)

\item[(ii)] A hypersurface has contact type as a coisotropic
submanifold if and only if it has contact type in the standard
sense. A hypersurface is stable as a coisotropic submanifold if and
only if it is stable as a hypersurface, i.e., there exists a vector
field $Z$ transverse to $M$ and such that
$\ker\left(\varphi_t^*\omega |_{M}\right)=\ker\left(\omega
|_{M}\right)$ for small $|t|$, where $\varphi_t$ is the flow of $Z$;
cf.\  \cite[p.\ 122]{HZ} and \cite{EKP}.
We will generalize this observation to coisotropic
submanifolds of codimension $k\geq 1$ in Proposition \ref{prop:stable}.

\item[(iii)] The product of stable submanifolds is also stable.  More
precisely, let $M_1\subset W_1$ and $M_2\subset W_2$ be stable. Then
$M_1\times M_2\subset W_1\times W_2$ is stable.  For instance, the
product of contact type hypersurfaces is a stable coisotropic
submanifold. The product $M\times S^1\subset W\times T^*S^1$ has
(restricted) contact type, provided that $M$ has (restricted) contact
type. However, unless $M_1$ or $M_2$ is one-dimensional, the product
$M_1\times M_2$ need not have contact type even if $M_1$ and $M_2$
have (restricted) contact type. This follows from Remark
\ref{rmk:cohomology} below.

\item[(iv)] A stable Lagrangian submanifold is necessarily
a torus as can be seen from Proposition \ref{prop:flow} or from Remark
\ref{rmk:cohomology}.
\end{enumerate}

\end{Example}

\begin{Example}
\label{ex:comm}
Let $M$ be a regular level set of the map $(K_1,\ldots,K_{k})\colon
W\to\R^{k}$ whose components are proper, Poisson--commuting
Hamiltonians. Assume furthermore that the Hamiltonian flows of
$K_1,\ldots,K_{k}$ generate an action of a torus $\T^k$ on $M$. Then
$M$ is stable.  (To see this, define $\alpha_i$ on $\T^k$-orbits,
i.e., the leaves of $\FF$, by $\alpha_i(X_{K_j})=\delta_{ij}$, where
$X_{K_j}$ is the Hamiltonian vector field of $K_j$, and extend these
forms to $\T^k$-invariant one-forms on $M$ in an arbitrary way. Then
$T\FF\subset \ker d\alpha_i$ and clearly
\eqref{eq:ct} is also satisfied.) Note that here, similarly to the
assertion of the Arnold--Liouville theorem, we require the
Hamiltonians to generate a torus action only on $M$, but not on the
entire ambient space $W$.
\end{Example}

\begin{Remark}
\label{rmk:cohomology}
The requirements that $M$ is stable or has contact
type impose severe restrictions on the topology of $M$ and $\FF$.
Namely, assume first that $M$ is stable.  Let $V$ be the vector space
formed by linear combinations $\alpha=a_1\alpha_1+\cdots +
a_k\alpha_k$, where $a_i$ are constants.  The forms $\alpha$ are
closed along $\FF$ and the natural map from $V$ to the foliated de
Rham cohomology $H^1_{dR}(\FF)$ along $\FF$, sending $\alpha$ to its
cohomology class, is a monomorphism. (Here, as above, $M$ is closed.)
In particular, $\dim
H^1_{dR}(\FF)\geq k$.  Indeed, $\alpha|_{\FF}=df|_{\FF}$ would imply
that $\alpha_x|_{T_x\FF}=0$ at a critical point $x$ of $f$, which in
turn means that $a_1=\ldots=a_k=0$ since the forms
$\alpha_1,\ldots,\alpha_k$ are, by \eqref{eq:ct}, linearly independent
in $T_x^*\FF$, cf.\ \cite[Remark 3]{Bo2}.

When $M$ has contact type, consider the vector space $V_0\subset V$
 formed by $\alpha=a_1\alpha_1+\cdots + a_k\alpha_k$ with
$a_1+\cdots+a_k=0$. The forms $\alpha$ are closed on $M$ and the
natural map $V_0\to H^1(M;\R)$ is a monomorphism. (As a consequence,
$\dim H^1(M;\R)\geq k-1$.) The proof of this observation due to Bolle,
\cite{Bo2}, is similar to the argument for stable
manifolds above. 
\ref{prop:flow}.
\end{Remark}

Let, as above, $M$ be a closed, stable, coisotropic
submanifold. Consider the product $M\times\R^k$.  Let
$(p_1,\ldots,p_k)$ be coordinates on $\R^k$ and let us use the same
symbols $\omega_0$ and $\alpha_i$ for differential forms on $M$ and
for their pull-backs to $M\times\R^k$. Then the form
\begin{equation}
\label{eq:normal-form}
\omega=\omega_0+\sum_{i=1}^k d(p_i\alpha_i),
\end{equation}
is symplectic near $M=M\times\{0\}$ in $M\times\R^k$. The normal bundle
to $M$ in $W$ is trivial, for it is isomorphic to $T^*\FF$, and thus
can be identified with $M\times\R^k$. Then, as is immediately clear
(see \cite{Bo1,Bo2}) from the Weinstein symplectic neighborhood
theorem (see, e.g., \cite{We77}), the local normal form of $\omega$
near $M$ is given by \eqref{eq:normal-form}.  From now on, we identify
a neighborhood of $M$ in $W$ with a neighborhood of $M$ in
$T^*\FF=M\times\R^k$ equipped with the symplectic form \eqref{eq:normal-form}.

\begin{Proposition}
\label{prop:flow}
Let $M$ be a stable coisotropic submanifold.

\begin{enumerate}

\item[(i)] The leaf-wise metric $(\alpha_1)^2+\cdots+(\alpha_k)^2$
on $\FF$ is leaf-wise flat.

\item[(ii)] The Hamiltonian flow of $\rho=(p_1^2+\cdots+p_k^2)/2$ is
the leaf-wise geodesic flow of this metric.

\end{enumerate}

\end{Proposition}

This proposition is essentially proved in \cite{Bo2}. Here we just briefly
outline the argument for the sake of completeness.

\begin{proof}
To prove the first assertion, note that by \eqref{eq:ct} the forms
$\alpha_i$ are linearly independent and leaf-wise closed. Thus, locally
on every leaf, their primitives form a coordinate system in which the
metric is isometric to the standard metric on $\R^k$. Let
$Y_1,\ldots,Y_k$ be the coordinate vector fields for this coordinate
system, i.e., $\alpha_i(Y_j)=\delta_{ij}$.  To prove the second
assertion, we now simply observe from \eqref{eq:normal-form} that the
Hamiltonian vector field $X_\rho$ of $\rho$ is given by
\begin{equation}
\label{eq:spray}
X_\rho=\sum_{i=1}^k p_i Y_i,
\end{equation}
which is the geodesic spray of the metric.
\end{proof}

The leaf-wise metric $\rho$ has some other relevant properties that
do not hold for leaf-wise flat metrics in general.
For instance, the length spectrum of $\rho$ is nowhere dense;
see Lemma~\ref{lemma:nowheredense}.

The next proposition, also quite elementary, clarifies the nature
of the stability condition and generalizes Example \ref{exam:Bolle}(ii).

\begin{Proposition}
\label{prop:stable} A coisotropic submanifold $M$ is stable if and
only if there exists a tubular neighborhood $M\times U$ of
$M=M\times\{0\}$ in $W$, where $U\subset\R^k$ is a neighborhood of
the origin, such that the submanifolds $M_p=M\times\{p\}$ are coisotropic
for $p\in U$ and $\ker\omega_p=\ker\omega_0$, where
$\omega_p=\omega|_{M_p}$.
\end{Proposition}

\begin{proof}
The fact that stability implies the existence of such a neighborhood
is a consequence of the normal form
\eqref{eq:normal-form}. Conversely, set $\alpha_i = \left(i_{\p/\p
p_i}\omega\right)|_M$. Then, \eqref{eq:ct} follows immediately from
the fact that $\p/\p p_1,\ldots,\p/\p p_k$ are linearly
independent. Furthermore,
$$
d\alpha_i=\left.\left(L_{\frac{\p}{\p p_i}}\omega\right)\right|_M
=\left.\frac{d}{d t}\omega_{p(t)}\right|_{t=0},
$$
where $p(t)=(0,\ldots, 0, t, 0,\ldots, 0)$ with $t$ in the $i$th slot, and
hence $\ker d\alpha_i\supset \ker \omega_0$ as required.
\end{proof}

\subsection{Hamiltonian displacement of coisotropic submanifolds}
\label{sec:statement}
To state our results on Hamiltonian displacement of stable coisotropic
submanifolds, we need to impose some natural conditions on the ambient
symplectic manifold $(W,\omega)$. Namely, in what follows $W$ is
always assumed to be \emph{symplectically aspherical}, i.e.,
$\omega|_{\pi_2(W)}=0=c_1|_{\pi_2(W)}$.

Furthermore, we require $W$ to be closed or \emph{geometrically
bounded} and \emph{wide}. The condition that $W$ is geometrically
bounded means that $W$ admits a complete metric which is compatible with
$\omega$ in a rather weak sense and has  injectivity radius bounded away
from zero and sectional curvature bounded from above; see \cite{AL}
for the precise definition. A symplectic manifold is said to be wide
if it admits a proper Hamiltonian bounded from below whose
Hamiltonian flow has no non-trivial contractible periodic orbits of
period less than or equal to one; see \cite{gu:new}. Among wide
manifolds are manifolds convex at infinity (e.g., cotangent bundles
to closed manifolds and $\R^{2n}$), twisted cotangent bundles, and non-compact
covering spaces of closed manifolds. In fact, the author is not aware of
any example of a geometrically bounded, open manifold that is not wide.

The essence of these requirements is that the standard machinery of
Floer homology is applicable to symplectically aspherical,
geometrically bounded manifolds (see, e.g., \cite{cgk,gg:new} and
Section \ref{sec:non-deg} below). Furthermore, one of the main tools
utilized in this paper is the technique of action selectors.  This
technique, developed for closed and convex at infinity symplectically
aspherical manifolds in \cite{FS,schwarz}, has recently been extended
to geometrically bounded, wide manifolds by G\"urel,~\cite{gu:new}.

Recall also that the energy of a compactly supported,
time-dependent Hamiltonian $H\colon [0,\,1]\times W\to \R$ is defined
as
$$
\parallel H \parallel =\int_0^1 (\max H_t-\min H_t)\, dt,
$$
where $H_t=H(t,\cdot)$. The Hamiltonian diffeomorphism $\varphi_H$, i.e.,
the time-one Hamiltonian flow of $H$, is said to displace $M$ if
$\varphi_H(M)\cap M=\emptyset$. When such a map $\varphi_H$ exists, we call
$M$ displaceable. For instance, every compact subset of $\R^{2n}$ is
displaceable.

Now we are in a position to state the main result of the paper.

\begin{Theorem}
\label{thm:main}
Let $W$ be symplectically aspherical, and closed or wide and
geometrically bounded. Let $M$ be a closed, stable, coisotropic
submanifold of $W$.
\begin{enumerate}

\item[(i)] Then there exists a constant $\Delta>0$ such that $\parallel H
\parallel > \Delta$ for any compactly supported Hamiltonian $H$ with
$\varphi_H$ displacing $M$.

\item[(ii)] When $M$ is displaceable,
there exists a loop $\eta$ tangent to $\FF$, contractible in $W$,
and bounding a positive symplectic area $A(\eta)$.

\item[(iii)] Moreover, if $M$ is displaceable and has restricted contact type,
there exists a loop $\eta$ tangent to $\FF$, contractible in $W$,
and such that $0<A(\eta)\leq \|H\|$.
\end{enumerate}
\end{Theorem}

\begin{Remark}
Note that the loop $\eta$ from assertions (ii) and (iii) is
necessarily not contractible in the leaf containing it.  Moreover,
$\eta$ is not contractible in the class of loops tangent to
$\FF$. This follows immediately from the observation that the area
spectrum of $M$ has zero measure.  (By definition, the area spectrum
of $M$ is the set formed by symplectic areas bounded by loops in $M$
that are tangent to $\FF$ and contractible in $W$.)  The same holds
for the curve $\gamma$ from Theorem \ref{thm:comm} below.
\end{Remark}

Theorem \ref{thm:main} will be proved
in Section \ref{sec:prfs}.  For stable coisotropic submanifolds in
$\R^{2n}$ assertions (ii) and (iii) were established by Bolle,
\cite{Bo1,Bo2}.

When $M$ is Lagrangian, and hence necessarily a torus, the first
assertion is a particular case of the Lagrangian intersection property
discussed in Section \ref{sec:intro}; see, e.g.,
\cite{Ch:Cambr,Ch,F:Morse,F:grad,F:c-l,Fl,Gr,LS:lagr,Oh1,Oh:cambridge,
Oh:disj,We73,We74} for similar and more general results. When $M$ is a
hypersurface, assertion (i) is trivial. However, in this case, we
prove a sharper theorem concerning leaf-wise intersections and
complementing the results of \cite{Ba,Dr,EH,Ho90,Li,Mo}, cf.\ Example
\ref{ex:leafwise} and Question~\ref{quest:lw-intersections}.

\begin{Theorem}
\label{thm:leafwise} Let $M$ be a connected closed hypersurface of
restricted contact type, bounding a domain $U$ in a subcritical Stein
manifold $W$.  Then for a compactly supported Hamiltonian
diffeomorphism $\varphi_H\colon W\to W$ with $\|H\|<\chom(U)$, there
exists a leaf $F$ of the characteristic foliation on $M$ such that
$\varphi_H(F)\cap F\neq\emptyset$.
\end{Theorem}

We refer the reader to Section \ref{sec:homcap} for the definition of
the homological capacity $\chom$. Theorem \ref{thm:leafwise} will be
established in Section \ref{sec:leafwise}. For $W=\R^{2n}$, this
result was proved in \cite{Ho90}.

As stated, with the upper bound on $\|H\|$, Theorem \ref{thm:leafwise} does
not hold for hypersurfaces that do not have contact
type.  In Example \ref{ex:leafwise-fail}, we construct a
Hamiltonian flow $\varphi^t$ on $\R^{2n}$ and a sequence of
coisotropic submanifolds $M_i$, $C^0$-converging to $S^{2n-1}$,
such that $M_i$ and $\varphi^{t_i}(M_i)$ have no leaf-wise intersections
for some sequence of times $t_i\to 0+$.

As is pointed out in Section \ref{sec:intro}, assertion (ii) of
Theorem \ref{thm:main} fails when the requirement that $M$ is stable
is dropped.  Indeed, for hypersurfaces in $\R^{2n}$, (ii) implies the
existence of a closed characteristic,
while counterexamples to the Hamiltonian Seifert conjecture show
that in general closed characteristics need not exist; see, e.g.,
\cite{Ci,Gi95,gi:bayarea,gg:ex2,ke:example} and references therein.
Recall also that the condition that $M$ is displaceable is essential
in (ii): the Liouville class of the zero section of a cotangent bundle
is zero.

Assertion (ii) gives an affirmative answer to Question
\ref{quest:Liouville} for stable coisotropic manifolds. Namely, assume
that $\omega$ is exact, i.e., $\omega=d\lambda$. Then the restriction
$\lambda|_{M}$ is closed along $\FF$, and hence the foliated Liouville
class $[\lambda|_{\FF}]\in H^1_{dR}(\FF)$ is defined.

\begin{Corollary}
\label{cor:liouville}
Assume in addition to the hypotheses of Theorem \ref{thm:main} that
$\omega$ is exact and $M$ is displaceable. Then $[\lambda|_{\FF}]\neq 0$.
\end{Corollary}

\subsection{Commuting Hamiltonians}
Let, as in Section \ref{sec:dense} and Example \ref{ex:comm}, $M$ be a
regular level $\vec{K}^{-1}(0)$ of the map
$\vec{K}=(K_1,\ldots,K_{k})\colon W\to\R^{k}$ whose components are
Poisson--commuting Hamiltonians.  (Now, in contrast with Example
\ref{ex:comm}, we do not assume that $\vec{K}$ gives rise to
a torus action on $M$.) As above, $W$ is
required to be symplectically aspherical, and closed or wide and
geometrically bounded. In the latter case, we also require the map
$\vec{K}$ to be proper (on its image) near $M$ to insure that the
coisotropic manifolds $M_a= \vec{K}^{-1}(a)$ are compact and close to
$M$ when $a\in\R^k$ is near the origin. Then, as in Section
\ref{sec:dense}, consider loops on $M_a$ which are tangent to the
characteristic foliation, contractible in $W$, but not contractible in
the leaf of the foliation.  Such a loop can be thought of as an
analogue of a closed characteristic on a hypersurface; see
Section \ref{sec:dense}.  Hence, the existence of such a loop 
can be interpreted as a generalization of the
dense existence theorem of Hofer and Zehnder and of Struwe, \cite{HZ:V,HZ,
St}, to the moment map $\vec{K}$. In Section~\ref{sec:prfs} we prove

\begin{Theorem}
\label{thm:comm}
Assume that $M$ is displaceable. Then, for a dense set of regular values
$a\in\R^k$ near the origin, the level set $M_a$ carries a closed
curve $\gamma$,  which is contractible in $W$, tangent to the characteristic
foliation $\FF_a$ on $M_a$, and bounds positive symplectic area.
\end{Theorem}

\section{Filtered Floer homology}
\label{sec:non-deg}
In this section we recall the definition of filtered Floer homology
for geometrically bounded symplectically aspherical manifolds, set
conventions and notation used in this definition, and revisit the
construction of action selectors for the manifolds in question. Most of
the results mentioned here are either well known or established
elsewhere or can be proved by adapting standard arguments. For this
reason, the proofs are omitted or just very briefly outlined; however,
in each case detailed references are provided although not necessarily
to the original proofs. We refer the reader to Floer's
papers \cite{F:Morse,F:grad,F:c-l,F:witten,Fl}, to
\cite{BPS,cgk,fh,FHS,Oh:cambridge, SZ}, or to \cite{HZ,mdsa-book,Sa}
for introductory accounts of the construction of Floer homology in
this setting.

\subsection{Preliminaries: notation and conventions}
Let $(W^{2n},\omega)$ be a symplectically aspherical manifold.
Denote by $\Lambda W$ the space of smooth contractible
loops $\gamma\colon S^1\to W$ and consider a time-dependent Hamiltonian
$H\colon S^1\times W\to \R$, where $S^1=\R/\Z$. Setting
$H_t = H(t,\cdot)$ for $t\in S^1$, we define the action functional
$A_H\colon \Lambda W\to \R$ by
$$
A_H(\gamma)=A(\gamma)+\int_{S^1} H_t(\gamma(t))\,dt,
$$
where $A(\gamma)$ is the negative symplectic area bounded by $\gamma$, i.e.,
$$
A(\gamma)=-\int_z\omega,
$$
where $z\colon D^2\to W$ is such that $z|_{S^1}=\gamma$.

The least action principle asserts that the critical points of $A_H$
are exactly contractible one-periodic orbits of the time-dependent
Hamiltonian flow $\varphi_H^t$ of $H$, where the Hamiltonian vector
field $X_H$ of $H$ is defined by $i_{X_H}\omega=-dH$. We denote the
collection of such orbits by $\PP_H$ and let $\PP^{(a,\,b)}_H\subset
\PP_H$ stand for the collection of orbits with action in the interval
$(a,\,b)$.  The action spectrum $\CS(H)$ of $H$ is the set of critical
values of $A_H$. In other words, $\CS(H)=\{ A_H(\gamma)\mid \gamma\in
\PP_H\}$. This is a zero measure set; see, e.g., \cite{HZ,schwarz}.

In what follows we will always assume that $H$ is compactly supported
and set $\supp H=\bigcup_{t\in S^1}\supp H_t$.  In this case, $\CS(H)$
is compact and hence nowhere dense.

Let $J=J_t$ be a time-dependent almost complex structure on $W$. A Floer
anti-gradient trajectory $u$ is a map $u\colon \R\times S^1\to W$
satisfying the equation
\begin{equation}
\label{eq:floer}
\frac{\p u}{\p s}+ J_t(u) \frac{\p u}{\p t}=-\nabla H_t(u).
\end{equation}
Here the gradient is taken with respect to the time-dependent
Riemannian metric $\omega(\cdot,J_t\cdot)$. This metric gives rise to
an ($L^2$-) Riemannian metric on $\Lambda W$ and \eqref{eq:floer} can
formally be interpreted as the equation ${\p u}/{\p s}=-\nabla_{L^2}
A_H(u(s, \cdot))$. In other words, $u$ is a trajectory of the
$L^2$-anti-gradient flow of $A_H$ on $\Lambda W$. In what follows, we 
denote the curve $u(s,\cdot)\in \Lambda W$ by $u(s)$.

The energy of $u$ is defined as
\begin{equation}
\label{eq:energy}
E(u)=\int_{-\infty}^\infty \left\|\frac{\p u}{\p s}\right\|_{L^2(S^1)}^2\,ds
=\int_{-\infty}^\infty \int_{S^1}\left\|\frac{\p u}{\p t}-J\nabla H (u)
\right\|^2 \,dt\,ds.
\end{equation}
We say that $u$ is asymptotic to $x^\pm\in \PP_H$ as $s\to\pm \infty$
or connecting $x^-$ and $x^+$ if $\lim_{s\to\pm\infty} u(s)=x^\pm$ in
$\Lambda W$. More generally, $u$ is said to be partially asymptotic to
$x^\pm\in \PP_H$ at $\pm\infty$ if $u(s^{\pm}_k)\to x^\pm$ for some
sequences $s^{\pm}_k\to\pm\infty$. In this case
$$
A_H(x^-)-A_H(x^+)=E(u).
$$
We denote the space of Floer trajectories connecting $x^-$ and
$x^+$, with the topology of uniform $C^\infty$-convergence on compact
sets, by $\CM_H(x^-,x^+,J)$ or simply by $\CM_H(x^-,x^+)$ when the
role of $J$ is not essential, even though $\CM_H(x^-,x^+,J)$ depends on
$J$.  This space carries a natural $\R$-action $(\tau\cdot
u)(t,s)=u(t,s+\tau)$ and we denote by $\hat{\CM}_H(x^-,x^+,J)$ the
quotient $\CM_H(x^-,x^+,J)/\R$.

Recall that $\gamma\in \PP_H$ is said to be non-degenerate if
$d\varphi_H\colon T_{\gamma(0)}W\to T_{\gamma(0)}W$ does not have one
as an eigenvalue. In this case, the so-called Conley--Zehnder index
$\MUCZ(\gamma)\in\Z$ is defined; see, e.g., \cite{Sa,SZ}.  Here we
normalize $\MUCZ$ so that $\MUCZ(\gamma)=n$ when $\gamma$ is a
non-degenerate maximum of an autonomous Hamiltonian with a small
Hessian. Assume that all periodic orbits with actions in the interval
$[A_H(x^+),A_H(x^-)]$, including $x^\pm$, are non-degenerate.
Then, for a generic $J$, suitable transversality conditions are
satisfied and $\CM_H(x^-,x^+,J)$ is a smooth manifold of dimension
$\MUCZ(x^+)-\MUCZ(x^-)$; see, e.g., \cite{fh,SZ} and references
therein.

\subsection{Filtered Floer homology and homotopy}
The objective of this section is two-fold. In its first part we
briefly outline the construction of filtered Floer homology
following closely \cite{fh}; see also \cite{BPS,cgk,schwarz}. (Note
that in the case of open geometrically bounded manifolds the
necessary compactness property of the moduli spaces of connecting
trajectories is guaranteed by Sikorav's version of the Gromov
compactness theorem; see \cite{AL}.)  In the second part, we discuss
properties of monotonicity maps. Here, we depart slightly from the
setting of \cite{BPS,cgk,fh}, for we need to consider also
non-monotone homotopies and to account for possible
non-compactness of $W$. For this reason, some of the proofs, still
quite standard, are outlined below.

Throughout the discussion of the filtered Floer homology
$\HF^{(a,\,b)}(W)$, we assume, when $W$ is open, that all the
intervals are in the positive range of actions, i.e., $a>0$ for
any interval $(a,\,b)$. This condition can be relaxed in some
instances and replaced by the requirement that $(a,\,b)$ does not
contain zero.  (The latter is clearly necessary with the
definitions we adopt here, for $H$ is assumed to be compactly
supported and thus $H$ always has trivial degenerate periodic orbits if
$W$ is open.)

\subsubsection{Filtered Floer homology: definitions}
Let $H$ be a compactly supported Hamiltonian on $W$. Assume that all
contractible one-periodic orbits of $H$ are non-degenerate if $W$ is
closed or that all such orbits with positive action are
non-degenerate when $W$ is open. This is a generic condition.
Consider an interval $(a,\,b)$, with $a>0$ when $W$ is open, such that
$a$ and $b$ are outside $\CS(H)$. Then the collection
$\PP^{(a,\,b)}_H$ is finite.  Assume furthermore that $J$ is regular,
i.e., the necessary transversality conditions are satisfied for moduli
spaces of Floer trajectories connecting orbits from
$\PP^{(a,\,b)}_H$. This is again a generic property as can be readily
seen by applying the argument from \cite{fh,FHS,SZ}.

Let $\CF_k^{(a,\,b)}(H)$ be the vector space over $\Z_2$ generated by
$x\in \PP^{(a,\,b)}_H$ with $\MUCZ(x)=k$. Define
$$
\p \colon \CF_k^{(a,\,b)}(H)\to \CF_{k-1}^{(a,\,b)}(H)
$$
by
$$
\p x=\sum_y \#\big(\hat{\CM}_H(x,y,J)\big)\cdot y.
$$
Here the summation extends over all $y\in \PP^{(a,\,b)}_H$ with
$\MUCZ(y)=\MUCZ(x)-1$ and $\#\big(\hat{\CM}_H(x,y,J)\big)$ is the
number of points, modulo 2, in $\hat{\CM}_H(x,y,J)$. (Recall that in
this case $\hat{\CM}_H(x,y,J)$ is a finite set by the compactness
theorem.) Then, as is well known, $\p^2=0$.  The resulting complex
$\CF^{(a,\,b)}(H)$ is the filtered Floer complex for $(a,\,b)$. Its
homology $\HF^{(a,\,b)}(H)$ is called the filtered Floer homology.
This is essentially the standard definition of Floer homology with
critical points outside $(a,\,b)$ being ignored.  Then
$\HF(H):=\HF^{(-\infty,\infty)}(H)$ is the ordinary Floer homology
when $W$ is compact. (As is well-known, $\HF_*(H)=H_{*+n}(W;\Z_2)$.)
In general, $\HF^{(a,\,b)}(H)$ depends on the Hamiltonian $H$, but
not on $J$; see Section~\ref{sec:homotopy}.

Let $a<b<c$. Assume that all of the above assumptions are satisfied for
all three intervals $(a,\,c)$ and $(a,\,b)$ and $(b,\,c)$. Then clearly
$\CF^{(a,\,b)}(H)$ is a subcomplex of $\CF^{(a,\,c)}(H)$, and
$\CF^{(b,\,c)}(H)$ is naturally isomorphic to the quotient complex
$\CF^{(a,\,c)}(H)/\CF^{(a,\,b)}(H)$. As a result, we have the long exact
sequence
\begin{equation}
\label{eq:seq}
\ldots\to \HF^{(a,\,b)}(H)\to\HF^{(a,\,c)}(H)\to \HF^{(b,\,c)}(H)\to\ldots .
\end{equation}
We will refer to the first map in this sequence as the inclusion map and
to the second one as the quotient map and to the whole sequence as
the $(a,\,b,\,c)$ exact sequence for $H$.

In the construction of the action selector for open manifolds given
in Section~\ref{sec:sel-open}, it will be convenient to work with
filtered Floer homology for the interval $(0,\,b)$ even though
$0$ is necessarily a critical value of the action functional. This
homology is defined as
\begin{equation}
\label{eq:zero-infty}
\HF^{(0,\,b)}(H)=\varprojlim_{\eps\to 0+}\HF^{(\eps,\,b)}(H),
\end{equation}
where the inverse limit is taken with respect to the quotient maps
and $\eps\to 0+$ in the complement of $\CS(H)$.
It is clear that this definition is equivalent to the original one when
$W$ is closed and $0$ is not in $\CS(H)$.

\subsubsection{Homotopy}
\label{sec:homotopy}
Let us now examine the dependence of $\HF^{(a,\,b)}(H)$ on $H$.
Consider a homotopy $H^s$ of Hamiltonians from $H^0$ to $H^1$. By
definition, this is a family of Hamiltonians parametrized by $s\in
\R$, and such that $H^s\equiv H^0$ when $s$ is large negative and
$H^s\equiv H^1$ when $s$ is large positive. Furthermore, let $J^s$ be a
family of $t$-dependent almost complex structures such that again
$J^s\equiv J^0$ when $s\ll 0$  and $J^s\equiv J^1$ when $s\gg 0$.
(In what follows, both the family of Hamiltonians
$H^s$ and the pair of families $(H^s,J^s)$ will be referred to as a
homotopy, depending on the context, when no confusion can arise.)
For $x\in \PP^{(a_0,\,b_0)}_{H^0}$ and $y\in
\PP^{(a_1,\,b_1)}_{H^1}$ denote by $\CM_{H^s}(x,y,J^s)$ the space of
solutions of \eqref{eq:floer} with $H=H^s$ and $J=J^s$.

Next we need to address the regularity issue. When $W$ is closed,
regularity of a homotopy $(H^s,J^s)$ is understood in the standard
sense, i.e., that the standard transversality requirements are met by the
homotopy $(H^s,J^s)$; see \cite{fh,FHS,SZ}. (This is a generic property,
\cite{fh,FHS,SZ}.) If $W$ is open, these conditions are never
satisfied, for Hamiltonians are compactly supported.  In this case, we
say that $(H^s,J^s)$ is regular as long as the transversality
requirements are met along all homotopy trajectories connecting
periodic orbits with positive action.  (This is also a generic
property; the argument of \cite{fh,FHS,SZ}  applies
to this case.)

When the transversality conditions are satisfied, $\CM_{H^s}(x,y,J^s)$
is a smooth manifold of dimension $\MUCZ(x)-
\MUCZ(y)$. Moreover, $\CM_{H^s}(x,y,J^s)$ is a finite set when
$\MUCZ(x)=\MUCZ(y)$.  Define the homotopy map
$$
\Psi_{H^0H^1}\colon \CF^{(a_0,\,b_0)}(H^0)\to \CF^{(a_1,\,b_1)}(H^1)
$$
by
$$
\Psi_{H^0H^1}( x)=\sum_y \#\big(\CM_{H^s}(x,y,J^s)\big)\cdot y.
$$
Here the summation is over all orbits $y\in \PP^{(a_1,\,b_1)}_{H^1}$ with
$\MUCZ(y)=\MUCZ(x)$ and $\#\big(\CM_{H^s}(x,y,J^s)\big)$
is the number of points, modulo 2, in this moduli space.

The map $\Psi_{H^0H^1}$ depends on the entire homotopy $(H^s,J^s)$
and in general is \emph{not} a map of complexes. However,
$\Psi_{H^0H^1}$ becomes a homomorphism of complexes when
$(a_0,\,b_0) = (a_1,\,b_1)$ and the homotopy is monotone decreasing,
i.e., $\p_s H^s\leq 0$ point-wise. Moreover, the induced map in
homology is then independent of the homotopy -- as long as the
homotopies are decreasing -- and commutes with the maps from the
exact sequence \eqref{eq:seq}.  (The reader is referred to, e.g.,
\cite{BPS,cgk,fh,Sa,SZ,schwarz,vi:functors}, for the proofs of these
facts for both open and closed manifolds.)  There are
other instances when the same is true. This is the case, for
instance, when the location of the intervals $(a_0,\,b_0)$ and
$(a_1,\,b_1)$ is compatible with the growth of the Hamiltonians in
the homotopy. We now analyze this particular case in more detail,
for it is essential for the proof of Proposition
\ref{prop:orbits-nondeg} below.

\begin{Definition}
\label{def:c-bound}
A homotopy $H^s$ is said to be \emph{$C$-bounded}, $C\in\R$, if
\begin{equation}
\label{eq:c-bound}
\int_{-\infty}^{\infty}\int_{S^1}\max_{W} \p_s H^s_t \,dt\,ds \leq C.
\end{equation}
\end{Definition}

It is clear that a $C$-bounded homotopy is also $C'$-bounded for
any $C'\geq C$. In what follows we will always assume that $C\geq 0$.

\begin{Example}
\label{ex:c-bound}
Let us give some examples of $C$-bounded homotopies.
\begin{enumerate}
\item[(i)] A monotone decreasing homotopy is $0$-bounded.

\item[(ii)] Every homotopy is $C$-bounded if $C$ is large enough.

\item[(iii)] Let $H^s$ be a \emph{linear homotopy} from $H^0$ to $H^1$, i.e.,
$$
H^s=(1-f(s))H^0+f(s)H^1,
$$
where $f\colon \R\to [0,1]$ is a monotone increasing compactly
supported function equal to zero near $-\infty$ and equal to one near
$\infty$. Then $H^s$ is $C$-bounded for any $C\geq
\int_{S^1}\max_W(H^1-H^0)\,dt$.
\end{enumerate}
\end{Example}

The following three observations (H0)--(H2) show that the standard
properties of homotopy maps in ordinary Floer homology extend to the
maps induced by $C$-bounded homotopies in filtered Floer homology. The
proof of (H0)--(H2) will further clarify the essence of Definition
\ref{def:c-bound}.  Note also that in (H0)--(H2) we will assume that
the end points of the intervals in question are outside the action
spectra of the Hamiltonians and, as above, $a>0$ when $W$ is open.
Furthermore, we will require all Hamiltonians to be non-degenerate (in
the positive action range for open manifolds) and the homotopies to be
regular.  However, this latter requirement is not essential, as we
will show in Section \ref{sec:inv}, and is not met in
Examples~\ref{ex:incl-quot} and~\ref{exam:incr}.

\begin{enumerate}
\item[(H0)] Let $H^s$ be a $C$-bounded homotopy. Then
$$
\Psi_{H^0H^1}\colon \CF^{(a,\,b)}(H^0)\to \CF^{(a+C,\,b+C)}(H^1)
$$
is a homomorphism of complexes for any interval $(a,\,b)$.
Hence, $\Psi_{H^0H^1}$ induces
a map in Floer homology, which we, abusing notation, also denote by
$\Psi_{H^0H^1}$. This map sends the $(a,\,b,\,c)$ exact sequence for $H^0$
to the $(a+C,\,b+C,\,c+C)$ exact sequence for $H^1$, i.e., on the level
of homology $\Psi_{H^0H^1}$ commutes with all maps in the long exact sequence
\eqref{eq:seq}.
\end{enumerate}

Here, and in (H1) and (H2) below, we suppress the dependence of this
map on the homotopy $J^s$ of complex structures. In fact, as we will
see in Section \ref{sec:inv}, the induced map on the level of homology is
independent of $J^s$.

\begin{Example}
\label{ex:incl-quot}
For any $C\geq 0$ the identity homotopy from $H$ to $H$
induces a map which is the composition
$$
\HF^{(a,\,b)}(H)\to \HF^{(a,\,b+C)}(H)\to \HF^{(a+C,\,b+C)}(H)
$$
of the inclusion and quotient maps from the exact sequences \eqref{eq:seq}
for the intervals $a<b<b+C$ and $a<a+C<b+C$, respectively. We will refer
to this map as the \emph{inclusion-quotient} map.

In view of Example \ref{ex:c-bound}(ii), note also that
$\CF^{(a+C,\,b+C)}(H)=0$ whenever $a+C>\max \CS(H)$, and hence the
assertion (H0) becomes trivial if, for a fixed homotopy, $C$ is
taken to be large enough. For instance, the inclusion-quotient map is
automatically zero if $a+C>b$.
\end{Example}

\begin{enumerate}
\item[(H1)] For a homotopy of $C$-bounded homotopies  $H^{s,\lambda}$,
i.e., a family of $C$-bounded homotopies parametrized by $\lambda\in [0,1]$,
the induced map
$$
\Psi_{H^0H^1}\colon \HF^{(a,\,b)}(H^0)\to \HF^{(a+C,\,b+C)}(H^1)
$$
is independent of the choice of homotopy in the family $H^{s,\lambda}$.
\end{enumerate}

Consider now a $C$-bounded homotopy $H^s$ from $H^0$ to $H^1$ and a
$C'$-bounded homotopy $G^s$ from $G^0=H^1$ to $G^1$. Define the gluing
or composition of these homotopies $H^s\#_R G^s$, where $R>0$ is large
enough, as $H^s\#_R G^s=H^{s+R}$ if $s\leq 0$ and $H^s\#_R
G^s=G^{s-R}$ if $s\geq 0$. (In other words, in $H^s\#_R G^s$ we
first, for $s\leq 0$, perform the homotopy obtained from $H^s$ by
shifting it sufficiently far to the left in $s$ and then, for $s\geq
0$, the homotopy resulting from $G^s$ by shifting it to the
right.) It is clear that $H^s\#_R G^s$ is $(C+C')$-bounded.

\begin{enumerate}
\item[(H2)] When $R>0$ is large enough, the map
$$
\Psi_{H^0G^1}\colon \HF^{(a,\,b)}(H^0)\to \HF^{(a+C+C',\,b+C+C')}(G^1)
$$
induced by the homotopy $H^s\#_R G^s$ is equal to the composition
$\Psi_{G^0G^1}\circ\Psi_{H^0H^1}$.
\end{enumerate}

\begin{Example}
\label{exam:incr}
Let $H^s$ be an increasing linear homotopy from $H^0$ to
$H^1\geq H^0$ and let $G^s$ be an arbitrary decreasing
homotopy from $G^0=H^1$ to $G^1=H^0$. By Example~\ref{ex:c-bound},
$H^s$ is $C$-bounded with $C=\int_{S^1}\max_W(H^1-H^0)\,dt$ and $G^s$
is 0-bounded. Hence, the composition of the homotopies is also
$C$-bounded and homotopic to the identity homotopy in the class of
$C$-bounded homotopies. By (H2) and (H3), this composition induces the
same map as the identity homotopy, i.e., the inclusion-quotient map from
Example \ref{ex:incl-quot}.
\end{Example}

\begin{proof}[Outline of the proof of (H0)--(H2)]
Consider first a solution $u$ of the homotopy Floer equation
\eqref{eq:floer} for a homotopy $H^s$, which is (partially)
asymptotic to $y$ at $\infty$ and $x$ at $-\infty$.
Then, as a direct calculation shows (see, e.g.,
\cite{schwarz}), we have
$$
A_{H^{s^+}}(u(s^+)) - A_{H^{s^-}}(u(s^-)) \leq
\int_{s^-}^{s^+}\int_{S^1} (\p_s H^s_t)(u(s))\,dt\,ds
-\int_{s^-}^{s^+}\int_{S^1} \left\|\p_s u\right\|^2\,dt\,ds.
$$
Hence,
\begin{equation}
\label{eq:dif-act-strong}
A_{H^1}(y)-A_{H^0}(x)\leq
\int_{-\infty}^{\infty}\int_{S^1}\max_W\p_s H^s_t\,dt\,ds - E(u),
\end{equation}
where $E(u)$ is still defined by \eqref{eq:energy}. In particular,
when $H^s$ is $C$-bounded, we infer that
\begin{equation}
\label{eq:dif-act}
A_{H^1}(y)-A_{H^0}(x)\leq C,
\end{equation}
by  Definition \ref{def:c-bound}, and
\begin{equation}
\label{eq:energy-action-bound}
E(u)\leq C + \left(A_{H^0}(x)-A_{H^1}(y)\right).
\end{equation}

Proving (H0), let us start with analyzing the case where $W$ is closed. Then,
by our assumptions, all periodic orbits of $H^0$ and $H^1$ are
non-degenerate and the homotopy is regular, and (H0) follows by the
standard gluing and compactness argument; see, e.g.,
\cite{Fl,fh,SZ,Sa,schwarz:book}. The only additional point to check is that
gluing and compactification of the moduli spaces in question do not
involve periodic orbits outside the range of actions. Assume, for
instance, that a homotopy connecting trajectory $u$ from $x\in
\PP^{(a,\,b)}_{H^0}$ to $z\in \PP^{(a+C,\,b+C)}_{H^1}$ with
$\MUCZ(z)=\MUCZ(x)+1$ is obtained by gluing an $H^0$-downward
trajectory from $x$ to $x'\in \PP^{(a,\,b)}_{H^0}$ with a homotopy
trajectory from $x'$ to $z$. Then, another component of the boundary of
$\CM_{H^s}(x,z,J^s)$ is a broken trajectory connecting first $x$ to
some $y\in \PP_{H^1}$ by a homotopy trajectory and then $y$ to $z$ by
an $H^1$-downward trajectory.  We need to check that $y\in
\PP^{(a+C,\,b+C)}_{H^1}$. It is clear that $A_{H^1}(y)\geq A_{H^1}(z)>
a+C$. Furthermore, by \eqref{eq:dif-act}, we have $A_{H^1}(y)\leq
A_{H^0}(x')+C<b+C$. The rest of the proof proceeds in a similar
fashion.

When $W$ is open, extra care is needed because the Hamiltonians
are compactly supported and hence there are degenerate critical
points. However, these points do not enter the calculation. Indeed,
let, for instance, $x$ and $z$ be as above. Note that the energy of
trajectories from $\CM:=\CM_{H^s}(x,z,J^s)$ is uniformly bounded due
to \eqref{eq:energy-action-bound}. By the compactness theorem, the
closure $\bar{\CM}$ of $\CM$ is compact in $C^{\infty}(\R\times
S^1,W)$ equipped with the topology of uniform
$C^{\infty}$-convergence on compact sets.  Then a trajectory $v\in
\bar{\CM}$ also has bounded energy and thus is partially asymptotic
to some orbits $x'\in \PP_{H^0}$ and $y'\in \PP_{H^1}$. (This can be
verified by adapting, for instance, the argument from the proof of
Proposition 10 on page 235 of \cite{HZ}.) Moreover, $A_{H^0}(x')\leq
A_{H^0}(x)$ and $A_{H^1}(y')\geq A_{H^1}(z)$; cf.\ \cite[p.\
66]{schwarz:book}.  (To prove, for example, the second of these
inequalities, pick $s_k\to\infty$ such that $v(s_k)\to y'$ and all
$s_k$ are large enough so that $H^{s_k}=H^1$. Then
$A_{H^1}(u(s_k))\geq A_{H^1}(z)$ for all $u\in\CM$. For every $k$ we
have a sequence $u_n\in \CM$ with $u_n(s_k)\to v(s_k)$. Hence,
$A_{H^{1}}(v(s_k))\geq A_{H^1}(z)$ and $A_{H^1}(y')=\lim
A_{H^{1}}(v(s_k))\geq A_{H^1}(z)$.) It follows immediately that $y'$
is a non-degenerate orbit, since $A_{H^1}(y')\geq A_{H^1}(z)>a+C>0$.
Moreover, $A_{H^0}(x')>0$, for otherwise we would have $a+C<
A_{H^1}(y')\leq C$ by \eqref{eq:dif-act}. Thus $x'$ is also a
non-degenerate orbit with positive action. Since $H^s$ is regular,
the standard description of $\bar{\CM}$ via broken trajectories
applies, and the proof is finished as for compact manifolds.

Properties (H1) and (H2) are established by a similar reasoning
invoking the standard gluing and compactification argument, checking
that the orbits are in the required ranges of action, and, when $W$
is open, verifying as above that the compactifications of the
relevant moduli spaces do not involve orbits with non-positive
action.
\end{proof}

\begin{Remark}
Assertions (H0)--(H2) still hold when the interval $(a+C,\,b+C)$ is
replaced by any interval $(a^1,\,b^1)$ with end points outside
$\CS(H^1)$ and such that $a^1\geq a+C$ and $b^1\geq b+C$. This
generalization, however, adds little new information, for then the
homotopy map $\HF^{(a,\,b)}(H^0)\to \HF^{(a^1,\,b^1)}(H^1)$ is
the composition of the homotopy $\HF^{(a,\,b)}(H^0)\to
\HF^{(a+C,\,b+C)}(H^1)$ with the inclusion-quotient map
$\HF^{(a+C,\,b+C)}(H^0)\to \HF^{(a^1,\,b^1)}(H^1)$.
\end{Remark}

\subsubsection{Invariance of filtered Floer homology}
\label{sec:inv}
Properties (H0)--(H2) of the homotopy maps have a number of standard
consequences (see, e.g., \cite{BPS,cgk,fh,Sa,SZ,schwarz:book,vi:functors}),
two of which, (H3) and (H4) below, are of particular relevance for this paper.

First we note that (H1) implies that on the level of homology the
map $\Psi_{H^0H^1}$ is independent of the homotopy of the almost
complex structures. Moreover, this map is actually independent of
the homotopy of the Hamiltonians as long as the latter is
$C$-bounded. (Indeed, any two $C$-bounded homotopies can be
connected by a linear in $\lambda$ family of $C$-bounded homotopies.
Furthermore, since regular homotopies are dense, a $C$-bounded
homotopy can be approximated by a regular $C'$-bounded homotopy with
$C'$ arbitrarily close $C$. Then, the assumption that the end-points
of the intervals are outside the action spectra guarantees that
the approximating homotopy still induces a map for the action ranges
$(a,\,b)$ and $(a+C,\,b+C)$, provided that $C'$ is close enough to
$C$.) Thus $\Psi_{H^0H^1}$ is in fact a canonical map depending only
on the Hamiltonians $H^0$ and $H^1$, the interval $(a,\,b)$ with
$a>0$ when $W$ is open, and the constant $C\geq 0$. This justifies
our neglect of the homotopy regularity condition in Examples
\ref{ex:incl-quot} and \ref{exam:incr}. Note, however, that this
canonical map is defined only when a $C$-bounded homotopy between
$H^0$ and $H^1$ exists which is not always the case. (For instance,
such a homotopy fails to exists with $C=0$ whenever $H^1>H^0$
point-wise.)

The next result, also quite standard, is the continuity property for
filtered homology; see \cite{BPS,cgk,fh,schwarz:book,vi:functors}.

\begin{enumerate}
\item[(H3)] Let $(a^s,\,b^s)$ be a family (smooth in $s$) of
non-empty intervals such that $a^s$ and $b^s$ are outside
$\CS(H^s)$ for some homotopy $H^s$ and such that $(a^s,\,b^s)$ is
independent of $s$ when $s$ is near $\pm\infty$ and, abusing
notation, equal to $(a^0,\,b^0)$ and $(a^1,\,b^1)$, respectively.
Then there exists an isomorphism of homology
\begin{equation}
\label{eq:isom}
\HF^{(a^0,\,b^0)}(H^0)\stackrel{\cong}{\longrightarrow} \HF^{(a^1,\,b^1)}(H^1).
\end{equation}
\end{enumerate}

We will outline the proof of (H3) below, after the next assertion (H4)
is stated.  When the interval is fixed and the homotopy is monotone
decreasing, the isomorphism \eqref{eq:isom} is in fact
$\Psi_{H^0H^1}$. We emphasize however that in general this is not the
case: the isomorphism is not induced by the homotopy $H^s$, for the
homotopy need not map $\CF^{(a^0,\,b^0)}(H^0)$ to
$\CF^{(a^1,\,b^1)}(H^1)$ even for a fixed interval of actions.  Leaving
aside the general question in what sense this isomorphism is canonical
and natural with respect to the long exact sequence \eqref{eq:seq}, we
focus on the particular case of a fixed interval.

\begin{enumerate}
\item[(H4)] Assume that $a<b<c$ are outside $\CS(H^s)$ for all
$s$. Then the isomorphisms \eqref{eq:isom} from (H3) for all three
intervals can be chosen to map the $(a,\,b,\,c)$ exact sequence \eqref{eq:seq} for
$H^0$ to the $(a,\,b,\,c)$ exact sequence for $H^1$.
\end{enumerate}

\begin{proof}[Outline of the proof of (H3) and (H4)]
Let us recall the construction of the isomorphism \eqref{eq:isom} for
a fixed interval, say, $(a,\,b)$. (Here we follow the argument from
\cite{BPS} essentially word-for-word.)  First, we break up $H^s$, up
to a homotopy of homotopies, into a composition of a finite sequence
of homotopies $K^s_i$ each of which is ``small''. By this we mean that
for all $i$ both $K_i^s$ and the inverse homotopy $K^{-s}_i$ are
$\eps$-bounded, where $\eps>0$ is so small that the intervals
$(a,\,a+2\eps)$ and $(b,\,b+2\eps)$ do not meet $\CS(K^s_i)$ for all
$s$.  Now it is sufficient to prove (H3) and (H4) for one homotopy
$K^s=K_i^s$. In this case, the isomorphism \eqref{eq:isom} is defined as
\begin{equation}
\label{eq:map}
\HF^{(a,\,b)}(K^0)\to \HF^{(a+\eps,\,b+\eps)}(K^1)
\stackrel{\cong}{\longrightarrow}
\HF^{(a,\,b)}(K^1),
\end{equation}
where the second arrow is the inverse of the quotient-inclusion map.
(The latter map is obviously an isomorphism due to the requirement imposed on
$\eps$.) To see that \eqref{eq:map} is an isomorphism,
observe that its inverse is given by
$$
\HF^{(a,\,b)}(K^1)
\stackrel{\cong}{\longrightarrow}
\HF^{(a+\eps,\,b+\eps)}(K^1)\to
\HF^{(a+2\eps,\,b+2\eps)}(K^0)
\stackrel{\cong}{\longrightarrow}
\HF^{(a,\,b)}(K^0),
$$
where the middle arrow is the map induced by the inverse homotopy
$K^{-s}$.  This concludes the proof of (H3) for a fixed interval.
The general case is established in a similar fashion.  Finally, from the
second part of (H0), it is clear that the maps \eqref{eq:map} defined
for all three intervals send the $(a,\,b,\,c)$ exact sequence for $K^0$
to that for $K^1$. This proves (H4).
\end{proof}

Another standard consequence of (H0)--(H2), along the lines of (H3) and
(H4), is that the filtered Floer homology can be defined ``by
continuity'' for any Hamiltonian $H$, not necessarily non-degenerate,
and any $a<b$ outside  $\CS(H)$. (See, e.g., \cite{vi:functors}.)
Moreover, the long exact sequence and (H0)--(H4) also hold in this
case.  Here, as above, we are assuming that $a>0$ when $W$ is open.

\subsection{The action selector}
\label{sec:selector}
In this section we recall the definition and the
properties of the action selector that are used in the proof of
Proposition \ref{prop:orbits-nondeg}. The constructions differ
somewhat depending on whether the manifold is open and wide or
closed. We start by dealing with the case of wide manifolds, closely
following \cite{gu:new}, for this case appears to be more relevant
to the displacement questions and, from the author's perspective,
more transparent. Then, we very briefly review the construction in
the case where $W$ is closed. Here we follow \cite{schwarz} with
some minor alterations. Note also that for $\R^{2n}$ and cotangent
bundles actions selectors were constructed in \cite{HZ} and
\cite{Vi:gen}, respectively, and the approach of \cite{schwarz} has
been extended to manifolds convex at infinity in \cite{FS}.

\subsubsection{The action selector for wide manifolds}
\label{sec:sel-open}
Assuming that $W^{2n}$ is geometrically bounded and wide (see Section
\ref{sec:statement} or \cite{gu:new}), let us recall from
\cite{gu:new} the definition of the action selector $\sigma(K)$ for a
compactly supported, non-negative Hamiltonian $K\colon S^1 \times W\to
\R$. It is easy to see that since $W$ is wide, there exists a smooth
compactly supported function $F\colon W\to [0,\,\infty)$ without
non-trivial contractible periodic orbits with period $T\leq 1$ and
such that $F\geq K$ point-wise. Without loss of generality, we may
assume that $\supp F$ is a smooth connected manifold with boundary and
that the restriction of $F$ to the interior of its support is a Morse
function with finitely many critical points.

Under these assumptions, we have
$$
\HF_*^{(0,\,\infty)}(F)\cong \HM_{*+n}^{(0,\,\infty)}(F)
\cong H_{*+n}(\supp F,\partial\supp F;\Z_2)
$$
and, in particular,
$$
\HF_n^{(0,\,\infty)}(F)\cong \Z_2.
$$
(Recall that the filtered Floer homology for $(0,\,\infty)$ is defined
by \eqref{eq:zero-infty}.)  Denote the generator of this group -- the
fundamental class -- by $[{\max}_F]$. A monotone decreasing homotopy
from $F$ to $K$ induces a map
$$
\Psi_{F,K}\colon \HF^{(0,\,\infty)}(F)\to\HF^{(0,\,\infty)}(K),
$$
independent of the choice of homotopy. (This follows from the results
of the previous section and \eqref{eq:zero-infty}.) Set
$$
[{\max}_K]=\Psi_{F,K}\left([{\max}_F]\right)\in\HF^{(0,\,\infty)}_n(K).
$$
It is not hard to show that $[{\max}_K]$ is well defined, i.e., independent
of the choice of $F$. Then, for $K\geq 0$, we set
$$
\sigma(K)=\inf\{a> 0\mid j^a_K\left([{\max}_K]\right)=0\}\in \CS(K),
$$
where
$$
j^a_K\colon \HF^{(0,\,\infty)}(K)\to \HF^{(a,\,\infty)}(K)
$$
is the quotient map.

Finally note that this definition makes sense even when $K$ is not
assumed to be non-negative. However, it is not clear whether the
resulting action selector has the required properties (S0)--(S4) stated
below. If $W$ is convex at infinity, the actions selector described
above coincides, for non-negative Hamiltonians, with the one of~\cite{FS}.

\subsubsection{The action selector for closed manifolds}
\label{sec:sel-closed}
Assume that $W^{2n}$ is closed. Then
for any Hamiltonian $K$, we have the canonical identification
$$
\HF_*(K)\cong H_{*+n}(W;\Z_2).
$$
Denoting by $[{\max}_K]$ the image of the fundamental class under
this isomorphism, set, following \cite{schwarz},
$$
\sigma(K)=\inf\{a\in\R\mid j^a_K\left([{\max}_K]\right)=0\}\in \CS(K),
$$
where
$$
j^a_K\colon \HF(K)\to \HF^{(a,\,\infty)}(K)
$$
is the quotient map. Observe that when $K\geq 0$ we can assume that
$a\geq 0$ by (S1) below, as in the definition of the action selector
for open manifolds.

Note that here, in contrast with \cite{schwarz}, we consider non-normalized
Hamiltonians, i.e., we are not assuming that $\int_W K_t\omega^n=0$. As
a consequence, neither $\CS(K)$ nor $\sigma(K)$ is uniquely determined by
$\varphi_K$.

Alternatively, one can define the class $[{\max}_K]$ as follows.
Let $F$ be a Morse function which is $C^2$-close to a constant.
Then, canonically,
$$
\HF_*(F)\cong\HM_{*+n}(F)\cong H_{*+n}(W;\Z_2).
$$
Let ${\max}_F\in{\operatorname{CM}}_{2n}(F)$ be the ``fundamental
class'', i.e., the sum of all local maxima. Then $[{\max}_F]$ is the
generator of $\HF_n (F)\cong \Z_2$. Now, for a fixed $K$, we take
$F\geq K$ and define $[{\max}_K]$ to be the image of $[{\max}_F]$
under the monotone decreasing homotopy map
$$
\Psi_{F,K}\colon \HF(F)\to\HF(K).
$$

\subsubsection{Properties of the action selector}
\label{sec:sel-prop}
The action selector $\sigma$, defined as above, has the following
properties for non-negative Hamiltonians, regardless of whether $W$
is wide or closed:

\begin{enumerate}

\item[(S0)] $\sigma$ is monotone, i.e., $\sigma(K)\leq\sigma(H)$,
whenever $0\leq K\leq H$ point-wise;

\item[(S1)] $0\leq\sigma(K)\leq E^+(K)$ for any $K\geq 0$,
where
$$
E^+(K)=\int_{S^1}\max_W K_t\,dt;
$$

\item[(S2)] $\sigma(K)>0$, provided that $K\geq 0$ is
not identically zero;

\item[(S3)] $\sigma(K)$ is continuous in $K$ in the $C^0$-topology;

\item[(S4)] $\sigma(K)\leq E^+(H)$, whenever $\varphi_H$ displaces
$\supp K$ and $H\geq 0$.
\end{enumerate}

We refer the reader to \cite{gu:new} for the proofs of (S1)--(S3) when
$W$ is wide and to \cite{schwarz} when $W$ is closed. As stated, (S4)
is established in \cite{FGS} for closed manifolds and in \cite{gu:new}
for wide manifolds in the stronger form $\sigma(K)\leq \sigma(H)$,
following the earlier versions of this upper bound from
\cite{FS,gi:alan,HZ,Oh:disj,schwarz,Vi:gen}.

Note also that if all contractible one-periodic orbits of $K$ are
non-degenerate, setting $A_K(c)=\sup_i A_K(c_i)$ for $c=\sum c_i\in
\CF^{(0,\,\infty)}(K)$, we have
$$
\sigma(K)=\inf_{[c]=[{\max}_K]} A_K(c).
$$
As a consequence, there exists a cycle $c$ (possibly non-unique)
with $A_K(c)=\sigma(K)$. We call such a cycle a \emph{carrier} of the
action selector.

\subsubsection{Homological capacity}
\label{sec:homcap}
Let $U\neq W$ be an open set in $W$, where $W$ is either closed or
geometrically bounded and wide. Recall that the homological capacity
$\chom(U)\in [0,\,\infty]$ of $U$ is defined as
$$
\chom(U)=\sup_{\supp K\subset U} \sigma(K).
$$
For a compact subset $Z$ of $W$ define the displacement energy
$e(Z)$ as $\inf \| H\|$, where the infimum is taken over all compactly
supported Hamiltonians $H$ with $\varphi_H(Z)\cap Z=\emptyset$, and
set $e(U)=\sup_{Z\subset U} e(Z)$.  Then $\chom(U)\leq e(U)$.  This
well-known fact (see, e.g., \cite{FS,gu:new,HZ,schwarz,Vi:gen}) follows, for
instance, from (S4). 

For a closed set $M\subset U$ (e.g., a closed submanifold), we set
$$
\chom(M)=\inf_{M\subset U} \chom(U).
$$
Note that $\chom(M)=0$ whenever $M$ is infinitesimally
displaceable. On the other hand, drawing on the results of
\cite{vi:functors}, one can expect that a closed Lagrangian
submanifold necessarily has positive homological capacity.
In Remark \ref{rmk:Dr}, we show that
this is also true for stable coisotropic submanifolds.  It seems to be
unknown whether the homological capacity of every closed hypersurface
in $\R^{2n}$ is positive.

\section{Displacement and connecting trajectories}
\labell{sec:nondeg}

Let, as above, $W$ be symplectically aspherical and either closed
or geometrically
bounded and wide. The key to the
proof of Theorem \ref{thm:main} is the following result.

\begin{Proposition}
\label{prop:orbits-nondeg}
Let $U$ be an open subset of $W$ displaced by a Hamiltonian diffeomorphism
$\varphi_H$. Furthermore, let $K\geq 0$ be a Hamiltonian supported in $U$ and
such that $\max K$ is large enough (e.g., $\max K>
\parallel H \parallel$ and $\max K> \max \CS(H)$, provided that $H\geq
0$). Assume furthermore that
\begin {itemize}

\item all one-periodic orbits from $\PP_K^{(0,\,\infty)}$ are non-degenerate;

\item for all $t\in [0,\,1]$ the functions $K_t$ attain
their maximum at the same point $p\in U$, which is thus a non-degenerate
maximum of $K_t$ for all $t$;

\item $K$ is autonomous (i.e., independent of $t$) near $p$.

\item the eigenvalues of the Hessian $d^2 K_t$ of $K_t$ at $p$ are
close to zero, and thus $K$ has no one-periodic orbits other than $p$
in a neighborhood of $p$.

\end {itemize}
Then the flow of $K$ has a contractible one-periodic orbit
$\gamma$ in $U$ which is connected to $p$ by a Floer
anti-gradient trajectory and such that
\begin{equation}
\label{eq:action}
0<A_K(\gamma)- \max K \leq \parallel H \parallel.
\end{equation}
\end{Proposition}

\begin{Remark}
Here we assume that $W$ is equipped with an almost complex structure,
possibly time-dependent, which is compatible with $\omega$ and such
that the regularity requirements are satisfied. Hence the Floer
complex of $K$ is defined.
\end{Remark}

\begin{Remark}
As stated, with the upper bound in \eqref{eq:action}, the proposition
is apparently new.  (This upper bound is absolutely crucial for the
proof of Theorem \ref{thm:main} and without it the proposition is
quite straightforward.)  However, it should be noted that in the
context of Hofer's geometry a number of somewhat similar results have
been established, often under much less restrictive conditions on the
ambient manifold; see \cite{En,KL,LM,MDS,Oh:chain}.  One can expect
that along the lines of some of these results the difference
$A_K(\gamma)-\max K$ can hypothetically be bounded from above in terms
of the Hofer norm of the time-one flow $\varphi_K$.
\end{Remark}

\emph{Outline of the proof.} Without loss of generality, we may assume
that $E^+(H)=\| H\|$. Since $\sigma(K)\leq \| H\|<K(p)=A_K(p)$, the
maximum $p$ is not homologically essential, i.e., there exists a
carrier $y$ of the action selector, not containing $p$, and the
``maximum cycle'' $z\in \CF^{(0,\,\infty)}_n(K)$ containing $p$ such
that $z-y=\p x$ for a chain $x\in \CF^{(0,\,\infty)}_{n+1}(K)$.  As a
consequence, there is a Floer connecting trajectory $u$ from some
orbit $\gamma$ in $x$ to $p$.  We can choose $x$ so that the class
$[x^c]$ induced by $x$ in $\HF^{(c,\,\infty)}_{n+1}(K)$ is non-zero
as long as $K(p)<c<A_K(\gamma)$ and hence for $K(p)<a<b<A_K(\gamma)$ the
quotient map $\HF^{(a,\,\infty)}(K)\to \HF^{(b,\,\infty)}(K)$ sending
$[x^a]$ to $[x^b]$ is also non-zero.

Note that, since $\varphi_H$ displaces $\supp K$, we have
$\HF^{(b,\,\infty)}(K\#H)= \HF^{(b,\,\infty)}(H)=0$ once
$b>\max\CS(H)$. Finally, arguing by contradiction, assume that
 the gap
between $K(p)$ and $A_K(\gamma)$ is greater than $E^+(H)$, i.e.,
$E(u)=A_K(\gamma)-K(p)> E^+(H)$.  Now we utilize the fact that the
homology of $K$ cannot be completely destroyed by a relatively small
perturbation $H$. This implies that $\HF^{(b,\,\infty)}(K\#H)\neq 0$
as along as $[b-E^+(H),\, b]\subset(\max
K,\, A_K(\gamma))$. Hence, when $\max K$ is large enough, the gap
cannot exceed $E^+(H)$ and $E(u)=A_K(\gamma)-K(p)> E^+(H)$.

\begin{proof}[Proof of Proposition \ref{prop:orbits-nondeg}]
First note that without loss of generality we may assume that $H$
meets the following requirements:
\begin{itemize}

\item $H$ is periodic in time and $H_0=H_1\equiv 0$;

\item $H$ is compactly supported and $\min_W H_t=0$ for all
$t\in S^1$.

\end{itemize}

The first requirement can be satisfied by reparametrizing $H$,
i.e., replacing $H$ by a new Hamiltonian of the form
$\lambda'(t)H(\lambda(t),x)$ where $\lambda\colon [0,\,1]\to [0,\,1]$
is a monotone increasing function identically equal to zero for $t$
near zero and to one for $t$ near one.  To have the second requirement
met, we replace $H$ by a Hamiltonian of the form $f_t\cdot
(H_t-\min_W H_t)$, where $f_t$ is a non-negative cut-off function
equal to one on $\varphi^t_H(U)$.  Note that neither of these
alterations changes $\| H\|$.

Both of the requirements are purely technical rather than essential.
The first one is needed to ensure that the composition $K\# H$ is
one-periodic in time. The second one enables us to treat
wide and closed manifolds in the same way, with only superficial
discrepancies.  In what follows we will assume that $H$ meets these
conditions. Note that since $\min_W H_t=0$, we have $H\geq 0$ and
$\| H\|=E^+(H)$.

The proof of the theorem proceeds slightly differently depending on
whether $W$ is wide  or closed. Below, we first consider the
case of wide manifolds and then indicate modifications needed
when $W$ is closed.

\emph{Wide manifolds.}
Let $F\geq K$ be an autonomous Hamiltonian as in the definition of
the action selector. Without loss of generality, we may assume that $F$
has a unique local maximum which is also located at $p$ and that
$F(p)=K(p)$. Then under a monotone homotopy from $F$ to $K$ the
Floer cycle $p={\max}_F\in \CF_n^{(0,\,\infty)}(F)$ is mapped to a
Floer cycle $z=z_1 +\cdots + z_m$ of $K$ with $z_1=p$. We emphasize that the
point
$p$ does occur in $z$. The reason is that the trivial Floer
connecting trajectory from $p$ for $F$ to $p$ for $K$ is the only
connecting trajectory between these two trivial orbits by the standard
energy argument. (Here we use the assumption that $F(p)=K(p)$.)
Furthermore, by utilizing the condition that the eigenvalues of
$d^2K_t$ are close to zero, it is easy to arrange the homotopy so that
this connecting trajectory is non-degenerate, cf.\ \cite{FHS}.

Note also that since the homotopy is decreasing, we have
$$
F(p)=K(p)=A_{K}(p)\geq A_{K}(z_i)\text{ for all $i=2,\ldots, m$},
$$
and hence
$$
A_{K}(z)=K(p)=\max K.
$$

Clearly, $[z]=[{\max}_{K}]$, but $z$ cannot be the carrier of the action
selector for $K$, for
$$
A_{K}(p)=\max K > E^+(H)\geq \sigma (K),
$$
where the last inequality follows from (S4). (In other words, $p$ is
not homologically essential.) Let $y\in \CF^{(0,\,\infty)}_n(K)$ be a
carrier of the action selector $\sigma(K)$. Then
$[y]=[{\max}_{K}]=[z]$ and $z-y=\p x$ for some $x\in
\CF_{n+1}^{(0,\,\infty)}(K)$.  Furthermore, it is clear from this
chain of inequalities that $p$ does not enter the cycle $y$ and, as a
consequence, $x$ contains a periodic orbit $\gamma$ connected with $p$
by a Floer anti-gradient trajectory.  Then
$$
0<A_K(\gamma)-\max K.
$$

However, to ensure that \eqref{eq:action} holds in its entirety, we
need to impose an additional requirement on the chain $x$ and the
orbit $\gamma$ which are in general  not unique.

To this end, for a given $x$ denote by $x_1$ an orbit, occurring in
$x$, which is connected to $p$ and has the smallest possible action
$A_{K}(x_1)$ among all such orbits in $x$. Then we chose $x$ with
$A_{K}(x_1)$ attaining the smallest value for all $x$ with $z-y=\p x$.
We will show that for any $x$ and $\gamma=x_1$, meeting this action
minimization condition,
\begin{equation}
\labell{eq:action2}
A_{K}(x_1)\leq \max K + E^+(H),
\end{equation}
which implies \eqref{eq:action}.

First let us consider the standard quotient map
$$
j\colon \HF_{n+1}^{(a,\,\infty)}(K)\to\HF_{n+1}^{(b,\,\infty)}(K).
$$
We claim that, provided that $x$ and $x_1$ are action minimizing as above,
\begin{equation}
\label{eq:nonzero}
j\neq 0, \text{ when $\max K< a < b < A_K(x_1)$}.
\end{equation}
(In particular, both of the groups are non-zero.)
To see this, denote by $x^a$ and $x^b$ the images of the chain $x$ in
$\CF_{n+1}^{(a,\,\infty)}(K)$ and $\CF_{n+1}^{(b,\,\infty)}(K)$,
respectively.  By definition, $j(x^a)=x^b$ on the level of
complexes. Furthermore, both $x^a$ and $x^b$ are closed, since $\p
x=z-y\in \CF^{(0,\,\,a)}_n(K)$.  Hence, it suffices to show that
$[x^b]\neq 0$ in $\HF_{n+1}^{(b,\,\infty)}(K)$.

Assume the contrary: there exists $w\in \CF_{n+2}^{(b,\infty)}(K)$ with
$\p w= x^b$. In other words, there exists a chain in
$w'\in\CF_{n+2}^{(0,\,\infty)}(K)$ such that
$x'=x-\p w'\in \CF_{n+2}^{(0,b)}(K)$. Then clearly $\p x'=z-y$ and
no orbit entering $x'$ has action in the interval $(b,\,\infty)$, where
$b< A_K(x_1)$, which is impossible due to our choice of $x$. This
contradiction completes the proof of \eqref{eq:nonzero}.

Proceeding with the proof of \eqref{eq:action2}, we again make use of
the fact that $U$ is displaced by $\varphi_H$. Denote by $K\# H$ the
Hamiltonian
$$
(K\# H)_t=K_t+H_t\circ(\varphi^t_K)^{-1}
$$
generating the time-dependent flow
$\varphi^t_K\circ\varphi^t_H$. Since $\varphi_H$ displaces $\supp K$,
the one-periodic orbits of $K\# H$ are exactly the one-periodic orbits
of $H$ and moreover $\CS(K\# H)=\CS(H)$, as is well known (see, e.g.,
\cite{HZ}). Thus, by the continuity property (H3) of Floer homology,
$$
\HF^{(b,\,b')}(K\# H)=\HF^{(b,\,b')}(H)
$$
for any $b<b'$ which are not in $\CS(H)$. In particular,
\begin{equation}
\label{eq:zero}
\HF^{(b,\,\infty)}(K\# H)\cong \HF^{(b,\,\infty)}(H)=0,
\text{ whenever $b>\max \CS(H)$.}
\end{equation}

Consider a linear monotone increasing homotopy from $K$ to $K\#H$.
(Recall that $H\geq 0$.)  By Example \ref{exam:incr}, such a homotopy
induces a map
$$
\Psi_{K,K\#H}\colon \HF^{(a,\,\infty)}(K)\to \HF^{(b,\,\infty)}(K\# H),
 \text{ where $b=a+ E^+(H) $,}
$$
whenever $a\not\in\CS(K)$ and $b\not\in\CS(K\# H)$. The
composition of this map with the map
$$
\Psi_{K\#H,K}\colon \HF^{(b,\,\infty)}(K\#H)\to \HF^{(b,\,\infty)}(K)
$$
induced by a monotone decreasing homotopy from $K\#H$ to $K$ is the
quotient map $j$. (See Example \ref{exam:incr}.)

To finish the proof in the case where $W$ is wide, assume that
\eqref{eq:action2} fails: $E^+(H) < A_K(x_1)-\max K$. Pick $a> \max K$
such that $b=a+ E^+(H) < A_K(x_1)$, and $a$ and $b$ are not in
$\CS(K)$ and also $b\not\in\CS(K\#H)$.  Then $j=
\Psi_{K\#H,K}\circ\Psi_{K,K\#H}$ factors through the group
$\HF^{(b,\,\infty)}(K\# H)$. Since $b>a>\max K> \max \CS(H)$, this
group is zero by \eqref{eq:zero}. Hence, $j=0$, which contradicts
\eqref{eq:nonzero}.

\emph{Closed manifolds}. Only two points of the above argument require
modifications when $W$ is closed.

The first of these is the definition of the cycle $z$. Now we take as
$F\geq K$ a function which is $C^2$-close to a constant and has a
unique local maximum, equal to $K(p)$, that is also located at
$p$. Then as above $z$ is the image of the Floer cycle $p={\max}_F\in
\CF_n^{(0,\,\infty)}(F)$ under a monotone homotopy from $F$ to $K$. It
is clear that $p$ occurs in $z$, i.e., $z=p+z_2 +\cdots + z_m$ for the
same reason as for open manifolds.

The second point of the proof that is not obvious when $W$ is closed,
for we work with non-normalized Hamiltonians, is the equality of
action spectra $\CS(K\#H)=\CS(H)$. To see that this is the case, let us
first recall how the action spectrum depends on the Hamiltonian. Let
$G^s$ be a homotopy of Hamiltonians such that the time-one flows
$\varphi_{G^s}$ are independent of $s$, i.e., the Hamiltonians $G^s$
determine the same element in the universal covering of the group
$\Ham(W)$. Then
\begin{equation}
\label{eq:spectra}
\CS(G^1)=\CS(G^0)+
\int_{S^1}\int_W \left(G^1_t-G^0_t\right)\omega^n\,dt.
\end{equation}
This can be established by arguing as in the proof of \cite[Lemma
3.3]{schwarz}.\footnote{The author is grateful to Felix Schlenk for
his help in clarifying the behavior of action spectra under
homotopy.}

Following, e.g., \cite{HZ} or \cite{FGS}, consider a homotopy from
$K=K^0$ to $K^1$ through reparametrizations $K^s_t=\lambda'(t)
K_{\lambda(t)}$ of $K$ such that $K^1_t\equiv 0$ when $t\in
[0,\,1/2]$. Let $H^s$ be a similar homotopy beginning with
$H^0=H$ and ending with $H^1$, and such that $H^1_t\equiv 0$ for $t\in
[1/2,\,1]$. Clearly, these homotopies do not change the action spectra.
Let $G^s=K^s\# H^s$. It is easy to see that
$\CS(K^1\#H^1)=\CS(H^1)=\CS(H)$.  (The first equality follows from the
fact that $K^1\#H^1$ and $H^1$ have literally the same one-periodic
orbits, with the same parametrizations, and the two Hamiltonians are
equal along these orbits.)  Furthermore, a direct calculation shows
that
\begin{eqnarray*}
\int_{S^1}\int_W G^1_t\omega^n\,dt
&=&
\int_{S^1}\int_W (K^1_t+H^1_t)\omega^n\,dt\\
&=&
\int_{S^1}\int_W (K_t+H_t)\omega^n\,dt\\
&=&
\int_{S^1}\int_W G^0_t\omega^n\,dt.
\end{eqnarray*}
Thus, by \eqref{eq:spectra}, $\CS(K^1\# H^1)=\CS(K\# H)$. As a consequence,
$\CS(H)=\CS(K\# H)$.

The rest of the proof proceeds exactly as in the case of open manifolds.
\end{proof}

\section{Connecting trajectories for degenerate displaceable Hamiltonians}
\label{sec:deg}
In this section we extend Proposition \ref{prop:orbits-nondeg} to
Hamiltonians whose one-periodic orbits are degenerate. This result will be
used in the proof of Theorem \ref{thm:main}.

\subsection{The space of finite energy anti-gradient trajectories}
Let, as above, $W$ be a geometrically bounded, symplectically
aspherical manifold and let $K$ be a compactly supported Hamiltonian
on $W$. Denote by $\CB=\CB(K)$ the space of contractible Floer
anti-gradient trajectories with finite energy, i.e., the space of
solutions of \eqref{eq:floer} such that $E(u)<\infty$ and $u(s)$ is
contractible.  We equip this space with the weak $C^\infty$-topology,
i.e., the topology of uniform $C^\infty$-convergence on compact
sets. Note that $\PP_K\subset\CB$.

Recall that by the compactness theorem any sequence of $u_i\in
\CB$ such that $u_i(0,0)$ is bounded contains a converging
subsequence.  In fact, $\CB$ is a locally compact, separable,
metrizable space, and the evaluation map
$$
\ev\colon \CB\to W,\text{ defined as $\ev(u)=u(0,0)$,}
$$
is continuous and proper; see, e.g., \cite{HZ} and references
therein.  Moreover, this map is a homeomorphism on a complement of a
compact subset of $\CB$.

For $u\in \CB$, set $A_K(u)$ to be the action of $K$ on the closed curve
$u(0)$. Then $A_K\colon \CB\to \R$ is a continuous function, which is
identically zero outside a compact set.

The anti-gradient flow $\Phi$ of $A_K$ on $\CB$ is defined as the shift in
the $\tau$-direction: $\Phi^\tau(u)(s)=u(s+\tau)$ for all $\tau\in\R$.
Obviously, the fixed points of this flow, i.e., finite energy
trajectories independent of $s$, are exactly the one-periodic orbits
of $K$. The function $A_K$ is decreasing along the orbits of
$\Phi^\tau$, i.e., $A_K(\Phi^\tau(u))\leq A_K(u)$ for any $\tau\geq
0$. Moreover, $A_K$ is strictly decreasing along non-trivial orbits:
$A_K(\Phi^\tau(u))< A_K(u)$ for any $\tau>0$, when $u$ is not a
fixed point of the flow.

Denote by $\Gamma(u)=\{\Phi^\tau(u)\mid \tau\in \R\}$ the orbit of
$\Phi$ through $u$. By definition, the limit set $\omega^+(u)$ is
formed by the limits of all converging sequences $\Phi^{\tau_i}(u)$
with $\tau_i\to\infty$. The limit set $\omega^-(u)$ is defined
similarly, but with $\tau_i\to-\infty$. The following properties of
the closure $\bG(u)=\Gamma(u)\cup\omega^+(u)\cup\omega^-(u)$ and of
the limit sets $\omega^\pm(u)$ are well known:
\begin{itemize}

\item the sets $\bG(u)$ and $\omega^\pm(u)$ are compact, connected, and
invariant under the flow $\Phi$;

\item the action functional $A_K$ is constant on
$\omega^\pm(u)$, and $A_K|_{\omega^+(u)}\equiv \min_{\bG(u)}A_K$
and $A_K|_{\omega^-(u)}\equiv \max_{\bG(u)}A_K$;

\item both of the limit sets $\omega^\pm(u)$ are non-empty and
entirely comprised of the fixed points of $\Phi$.

\end{itemize}
It is clear that $u$ is partially
asymptotic to $x^\pm$ as $s\to\pm\infty$, i.e.,
$u(s_i^\pm)\to x^\pm$
in $C^\infty(S^1,W)$ for some sequences $s_i^\pm\to \pm\infty$,
if and only if $x^\pm\in\omega^\pm(u)$. (Here we treat $x^\pm$ simultaneously
as periodic orbits (elements of $\PP_K$) and as finite energy trajectories
(elements of $\CB$).)  If
$x^\pm$ are non-degenerate, we necessarily have
$\omega^\pm(u)=\{x^\pm\}$. Otherwise, the limit sets can be quite
large and $x^\pm$ are not unique.

\subsection{Existence of connecting trajectories for degenerate Hamiltonians}
In this section we prove an analogue of Proposition
\ref{prop:orbits-nondeg} for degenerate Hamiltonians. Let, as in
Section \ref{sec:nondeg}, $W$ be a symplectically aspherical manifold, which
is either closed or geometrically bounded and wide.

\begin{Proposition}
\label{prop:orbits-deg}
Let $U$ be an open subset of $W$ displaced by a Hamiltonian
diffeomorphism $\varphi_H$. Furthermore, let $K\geq 0$ be a
Hamiltonian supported in $U$ and such that for all $t\in [0,\,1]$ the
functions $K_t$ attain their maxima at the same connected set
$M\subset U$ and are autonomous near $M$. Assume also that $K$
is $C^2$-close to a constant on a sufficiently small neighborhood of $M$. Then,
provided that $\max K$ is large enough, there exists $u\in\CB(K)$
partially asymptotic to a point of $M$ and to a contractible
one-periodic orbit $\gamma$ such that \eqref{eq:action} holds:
$$
0<A_K(\gamma)- \max K \leq \parallel H \parallel.
$$
\end{Proposition}

\begin{Remark}
\label{rmk:orbits-deg}
In the context of this proposition, $W$ is equipped with an arbitrary
almost complex structure $J$ compatible with $\omega$ in the sense of
geometrically bounded manifolds and independent of time at infinity.
Note also that \eqref{eq:action} guarantees that $A(\gamma)>0$ and,
in particular, $\gamma$ is a non-trivial one-periodic orbit when $K$
is autonomous. Furthermore, $A_K(\gamma)- \max K=E(u)$.
\end{Remark}

\begin{Remark}
Proposition \ref{prop:orbits-deg} will be applied in the setting where
$M$ is a Morse--Bott non-degenerate critical set of $K$. Then it might
be possible to conclude that $u$ is genuinely asymptotic to a point on
$M$ as $s\to\infty$; cf.\ \cite{poz}. However, the remaining
one-periodic orbits of $K$ are still very degenerate and in this case
$u$ need not be truly asymptotic to any orbit at $-\infty$. (This phenomenon,
occurring already for the anti-gradient flow of a function on a smooth
manifold, is overlooked in \cite{Bo2}. However, the argument of \cite{Bo2}
still goes through for partial asymptotics.)

Note also that some degree of control over periodic orbits of $K$ near
$M$ is essential. Here, it is achieved through requiring $K$ to be
$C^2$-close to a constant near $M$ and thus to have no one-periodic
orbits near $M$ with action greater than $\max K$. Without a
restriction on the behavior of $K$ near $M$ the proposition probably
fails.
\end{Remark}

The idea of the proof is, of course, to approximate $K$ by a sequence
$K_l\to K$ of non-degenerate Hamiltonians satisfying the hypotheses of
Proposition \ref{prop:orbits-nondeg} and to define the orbit $u$ as the
limit of a sequence of trajectories $u^0_l$ for $K_l$ such as in Proposition
\ref{prop:orbits-nondeg}. The nuance is that an arbitrary sequence
$u^0_l$ need not have an orbit with the required properties
as its limit point. (For instance, $u^0_l$ can converge to a point on $M$
or to a trajectory which is not partially asymptotic to a point on $M$.)
However, the trajectories $u^0_l$ are not unique and can
be replaced by $u_l=\Phi^{\tau_l}(u^0_l)$ for any sequence of shifts
$\tau_l$. We show that $u_l$ converges to the required $u$ for some
sequence $\tau_l$, using an elementary, point--set topological
argument.

\begin{proof}
Pick a sequence of non-degenerate Hamiltonians $K_l$ which are supported
in $U$ and $C^\infty$-converge to $K$ and satisfy the hypotheses of
Proposition \ref{prop:orbits-nondeg}. Without loss of generality, we may
assume that all $K_l$ attain their maxima at the same point $p\in M$ and
that $\max K_l=K_l(p)=\max K$.

Furthermore, we can approximate $J$ by almost complex structures $J_l$
(possibly time-dependent) which are compatible with $\omega$ and equal
to $J$ at infinity and such that the pairs $(K_l,J_l)$ are regular.
Let $u_l^0$ be an anti-gradient trajectory for $K_l$ whose
existence is established in Proposition \ref{prop:orbits-nondeg}:
$u_l^0$ connects a contractible one-periodic orbit $\gamma^0_l$
of $K_l$ to $p$ and \eqref{eq:action} holds. By passing if necessary
to a subsequence, we may assume that the sequence $\gamma_l^0$ converges to a
one-periodic orbit $\gamma^0$ of $K$.

Recall that $K$ is assumed to be autonomous and $C^2$-close to a
constant near $M$.  It follows that every one-periodic orbit of $K$,
with action in the interval $[\max K, \infty)$, that meets a sufficiently
small neighborhood $V$ of $M$ must be a point of $M$. Furthermore,
$K_l$ are also $C^2$-close to a constant on $V$. As a consequence,
every one-periodic orbit of $K_l$ entering $V$ is trivial and we may
require, in addition, that the actions of $K_l$ on these orbits (with
the exception of $p$) are strictly smaller than $\max K_l$. (It is
easy to show that the approximations $K_l$ with this property do
exist.)  Since $K_l(\gamma^0_l)>\max K_l$ by \eqref{eq:action}, we
conclude that $\gamma^0_l$ does not enter $V$. Hence, $\gamma^0$ also
lies entirely outside $V$.

Recall that by the compactness theorem any sequence $u_l\in\CB(K_l)$
contains a converging subsequence, provided that the sequence
$\ev(u_l)$ is bounded in $W$. (Here convergence is again understood as
$C^\infty$-convergence uniform on compact subsets of $\R\times S^1$,
i.e., in the weak $C^\infty$-topology on the space $\CC$ of smooth
maps $\R\times S^1\to W$.)

Set $\Gamma_l=\Gamma(u^0_l)$ and denote by $\Sigma\subset \CC$ the set of all
limit points of sequences $u_l\in \Gamma_l$ (or, equivalently,
$u_l\in\bG_l$). This set has the following properties:
\begin{enumerate}

\item[(a)] the set $\Sigma$ contains $p$ and $\gamma_0$;

\item[(b)] the sets $\bG_l$ converge to $\Sigma$ in the Hausdorff
topology: for every neighborhood $\CU$ of $\Sigma$ in $\CC$ we have
$\bG_l\subset\CU$ when $l$ is large enough;

\item[(c)] the set $\Sigma$ is connected, compact, and invariant under the
flow $\Phi$;

\item[(d)] the flow $\Phi$ on $\Sigma$ is non-trivial;

\item[(e)] the action functional $A_K$ is not constant on $\Sigma$, and
$\min_\Sigma A_K\mid =A_K(p)=\max K$
and $\max_\Sigma A_K\leq \max K+\parallel H\parallel$.

\end{enumerate}

\begin{proof}[Proof of (a)--(e)]
The first assertion (a) is obvious. To prove (b), we just note that
otherwise we would have a sequence $u_{l_i}\in\bG_{l_i}\ssminus\CU$
for some $l_i\to \infty$. By compactness, $u_{l_i}$ must have a limit
point outside  $\CU$ and hence not in $\Sigma$, which is
impossible. In (c), only the fact that $\Sigma$ is connected requires
a proof. This readily follows from (b) and from the fact that all $\bG_l$
contain the point $p$. Assertion (e) is obvious
except for the statement that $A_K$ is non-constant, which is a
consequence of (d).  To prove (d), observe that $\ev(\Sigma)$ contains
both $p\in M\subset V$ and $\gamma^0(0)\not \in V$.  By (c),
$\ev(\Sigma)$ is connected and therefore there exists $u\in\Sigma$
such that $u(0,0)=\ev(u)\in V\ssminus M$. Since no one-periodic orbit
of $K$ with action in $[\max K, \infty)$, other than the points of
$M$, enter $V$, we conclude that $u(0)$ is not a contractible
one-periodic orbit of $K$ and, in particular, the entire anti-gradient
trajectory $u$ cannot be a contractible one-periodic orbit of
$K$. Thus $u$ is not a fixed point of the flow $\Phi$.
\end{proof}

Let $\Sigma_{\min}$ be the set of $u\in\Sigma$ at which
$A_K|_{\Sigma}$ attains its minimum $\max K$. Note that
$\Sigma_{\min}$ is entirely comprised of fixed points of $\Phi$ or,
equivalently, of periodic orbits of $K$.  Regarding $M$ as a subset of
$\CB(K)$, set $C=M\cap \Sigma_{\min}=M\cap \Sigma$. Then $C$ is a
compact, proper subset of $\Sigma$. Indeed, compactness of $C$ is
obvious.  By (e), $\Sigma_{\min}\neq \Sigma$, and hence $C\neq
\Sigma$.  On the other hand, $C\neq\emptyset$, for $p\in C$. Next, we
claim that
\begin{enumerate}

\item[(f)] the set $C$ is a union of connected components of $\Sigma_{\min}$,
i.e., $\Sigma_{\min}\ssminus C$ is closed.

\end{enumerate}

\begin{proof}[Proof of (f)]
Indeed, assume the contrary.  Then, there exists a sequence of
periodic orbits $u_l\in \Sigma_{\min}\ssminus C$ converging to a point
of $C$. As a consequence, $\ev(u_l)$ converges to a point of $M$. Then
$u_l$ must be a trivial periodic orbit when $l$ is large enough.  (For
$u_l$ enters $V$.) In addition, $A_K(u_l)=\max K$ and thus $u_l$ is a
point of $M$. This is impossible since $u_l\in \Sigma_{\min}\ssminus C$.
\end{proof}

Set $\CN_\eps=\{u\in\Sigma\mid A_K(u)<\max K+\eps\}$. Fix a connected
component $C_0$ of $C$ and let $\CU_\eps$ be the connected component
of $\CN_\eps$ that contains $C_0$.  Clearly, the sets $\CU_\eps$ with
$\eps>0$ are open, nested, and invariant under the positive flow
$\Phi^{\tau\geq 0}$. Our next goal is to prove
\begin{enumerate}

\item[(g)] the open sets $\CU_\eps$ with $\eps>0$ form a fundamental
system of neighborhoods of $C_0$, i.e., for any open set $\CU\supset C_0$,
we have $\CU_\eps\subset \CU$, when $\eps>0$ is sufficiently small.

\end{enumerate}

\begin{proof}[Proof of (g)]
Let us first show that the sets $\CN_\eps$ form a fundamental system
of neighborhoods of $\Sigma_{\min}$ for $\eps>0$. Assume the
contrary. Then for some open set $\CW\supset \Sigma_{\min}$, there
exists a sequence $u_l\in \CN_{\eps_l}\ssminus \CW$ with $\eps_l\to
0+$.  Passing if necessary to a subsequence, we have $u=\lim u_l \in
\Sigma\ssminus \CW$ by the compactness theorem, and hence
$u\not\in\Sigma_{\min}$.  This is impossible, for $A_K(u)=\lim
A_K(u_l)=\min A_K|_{\Sigma}$ and thus $u\in \Sigma_{\min}$.

Let now $\CU$ be a neighborhood of $C_0$. We need to show
that $\CU_\eps\subset \CU$, when $\eps>0$ is small. Assume that this
is not the case: the sets $\CU_{\eps_i}$ are not entirely contained in
$\CU$ for some sequence $\eps_i\to 0+$.  (Since the family
$\CU_{\eps}$ is nested, this is true for all $\CU_{\eps}$, but we
prefer to work with a sequence.) Therefore, the intersection
$\bigcap_i \overline{\CU}_{\eps_i}\subset \Sigma_{\min}$ is connected
since $\CU_{\eps_i}$ are connected,
contains $C_0$, and a point of $\Sigma_{\min}\ssminus C$. This is impossible
due to (f).
\end{proof}

Now we are in a position to finish the proof of Proposition
\ref{prop:orbits-deg}. Utilizing (g), pick $\eps>0$ so small that
$\ev(\CU_\eps)\subset V$.  Then every fixed point of $\Phi$ in
$\CU_\eps$ must belong to $C_0$. (Otherwise, there would be a one-periodic
orbit with action in $[\max K, \infty)$, other than a point of $M$,
entering $V$.)  Thus $\omega^+(u)\in C_0$ for any $u\in
\CU_\eps$. Observe now that $\CU_\eps \ssminus C_0\neq \emptyset$, since
$\Sigma$ is connected and $C_0$ is closed.  Let $u\in \CU_\eps\ssminus
C_0$. Then $u$ is a non-trivial, anti-gradient trajectory partially
asymptotic to a point in $M$ at $\infty$. Therefore,
$$
\max K=\min_{\Sigma} A_K < A_K(u).
$$
As a consequence, for any $\gamma\in \omega^-(u)$ we have
$$
0<A_K(\gamma)-\max K\leq\parallel H\parallel,
$$
where the second inequality follows from (e).  This completes the
proof of \eqref{eq:action} and thus the proof of the proposition.
\end{proof}

Proposition \ref{prop:orbits-deg} has a counter-part asserting the
existence of homotopy connecting trajectories ``transferring''
action selectors; cf.\ \cite{cgk,ke:new}. This is a much more
standard result and we treat it in lesser detail.

\begin{Proposition}
\label{prop:orbits-deg2}
Let $U$ be an open subset of $W$ displaced by a Hamiltonian
diffeomorphism $\varphi_H$. Furthermore, let $K\geq 0$ be a
Hamiltonian supported in $U$ and
let $f\geq 0$ be a $C^2$-small autonomous
Hamiltonian such that $f\leq K$.  Consider a monotone decreasing homotopy
from $K$ to $f$. Then there exists a homotopy trajectory $u$ partially
asymptotic to a point $p\in\PP(f)$ at $\infty$ and to
$\gamma\in\PP(K)$ at $-\infty$ and such that
\begin{itemize}
\item $f$ attaints its maximum at $p$,
\item $A_K(\gamma)\leq\sigma(K)$ and $E(u)\leq A_K(\gamma)-f(p)\leq \| H\|$.
\end{itemize}
\end{Proposition}

\begin{Remark}
When the Hamiltonians $K$ and $f$ are non-degenerate and the homotopy
is regular, the inequality $A_K(\gamma)\leq\sigma(K)$ can be replaced
by equality $A_K(\gamma)=\sigma(K)$.
\end{Remark}

\begin{proof}[Outline of the proof]
First note that similarly to the proof of Proposition
\ref{prop:orbits-nondeg} we may assume without loss of generality that
$H$ is compactly supported and $\min H_t=0$. In particular, $H\geq 0$
and $\| H\|=E^+(H)$. Furthermore, the lower bound $E(u)\leq
A_K(\gamma)-f(p)$ follows immediately from
\eqref{eq:energy-action-bound}.

When $K$ and $f$ are non-degenerate and the homotopy is regular the
assertion is well known; see, e.g., the proof of Proposition
\ref{prop:orbits-nondeg} or
\cite{FS,gi:alan,gu:new,ke:new,KL,Oh:chain,schwarz,Vi:gen} to mention just a
few sources where similar results have been proved.  Moreover, in this
case we have $A_K(\gamma)=\sigma(K)$. Indeed, in the notation of the
proof of Proposition \ref{prop:orbits-nondeg}, $\gamma$ is an orbit in
the chain $y$, a carrier of the action selector, such that
$A_K(\gamma)=\sigma(K)$.

To deal with the general case, we argue as in the proof of Proposition
\ref{prop:orbits-deg} and approximate $f$ and $K$ by non-degenerate
Hamiltonians $f_l\to f$ and $K_l\to K$ meeting the requirements of
Proposition \ref{prop:orbits-deg2}. We also approximate the homotopy
$K^s$ from $K$ to $f$ by regular monotone decreasing homotopies
$K^s_l$ from $K_l$ to $f_l$. Furthermore, we may assume that $\max
f_l=\max f$.  Applying the non-degenerate case of the proposition to
$K_l$ and $f_l$, we obtain a critical point $p_l$, an orbit $\gamma_l$,
and a homotopy connecting trajectory $u_l$.

Since $E(u_l)\leq \| H\|$, the compactness theorem implies that the
sequence $u_l$ contains a converging subsequence. Passing to this
subsequence and taking the limit, we obtain a homotopy connecting
trajectory $u$ from $K$ to $f$ such that $E(u)\leq \| H\|$. Then, by
compactness again, $u$ is partially asymptotic to $p\in\PP_f$ at
$\infty$ and to $\gamma\in\PP_K$ at $-\infty$. Note that although we
can assume that $\gamma_l$ converges to an orbit of $K$ due to the
Arzela--Ascoli theorem, we cannot claim that $\gamma_l\to\gamma$ and
hence cannot conclude that $\sigma(K)=A_K(\gamma)$. However,
\begin{equation}
\label{eq:f}
f(p)=A_f(p)\geq \max f = f_l(p_l)
\end{equation}
and
\begin{equation}
\label{eq:K}
A_K(\gamma)\leq \lim_{l\to\infty} A_{K_l}(\gamma_l)=\sigma(K).
\end{equation}
To prove \eqref{eq:K}, we argue as in
the proof of (H0)--(H2) in Section \ref{sec:homotopy}. Namely,
\begin{eqnarray*}
A_K(\gamma)=\sup_s A_{K^s}(u(s))
&=&
\sup_s\lim_{l\to\infty} A_{K^s_l}(u_l(s))\\
&\leq&
\lim_{l\to\infty} \sup_s A_{K^s_l}(u_l(s))\\
&=&
\lim_{l\to\infty} A_{K_l}(\gamma_l)\\
&=&
\lim_{l\to\infty} \sigma(K_l)=\sigma(K),
\end{eqnarray*}
where the last equality follows from (S3), continuity of the action selector.
The first inequality, \eqref{eq:f}, is established in a similar fashion
and immediately implies that $f(p) =\max f$. Finally, from \eqref{eq:f} and
\eqref{eq:K}, we infer that
$$
A_K(\gamma)-f(p)\leq
\lim_{l\to\infty}\left(A_{K_l}(\gamma_l)-f_l(p_l)\right)
\leq \| H\|,
$$
which concludes the proof.
\end{proof}

\begin{Remark}
Propositions \ref{prop:orbits-deg} and \ref{prop:orbits-deg2} (and
their non-degenerate counterparts) represent two different Floer
homological (broadly understood) approaches to proving the existence
of (infinitely many) periodic orbits of Hamiltonians, either in the
autonomous case (the Weinstein conjecture and the almost existence
theorem) or for time-dependent Hamiltonians (the Conley
conjecture). One approach comprises a class of methods that lead to
``low-lying'' orbits with action smaller than the displacement energy
as in Proposition \ref{prop:orbits-deg2}. These methods are utilized
in, for instance,
\cite{cgk,FS,fhw,FGS,gi:alan,gu:new,HZ,Vi:gen,vi:functors}. The second
class of methods detects orbits lying above the action selector value
as in Propositions \ref{prop:orbits-nondeg} and \ref{prop:orbits-deg}.
In this class are the variational methods of \cite{HZ:cap}, the Floer
homological results of \cite{ke:new,gg:new} and apparently some of the
results utilizing Hofer's geometry, e.g., \cite{Oh:chain,schl}.  It
goes without saying that many methods do not fit into this crude
classification. These include, for instance, the equivariant (Floer)
homological methods (see, e.g., \cite{HZ:V,vi,vi:jdg}) and holomorphic
curve methods, e.g., \cite{HV:spheres,Lu,Lu2,LiTi}.
\end{Remark}

\section{Displacement energy for stable manifolds}
\label{sec:prfs}
Let, as in Section \ref{sec:statements}, $M$ be a stable, closed, coisotropic
submanifold of $W$. Recall from Section \ref{sec:ct}
that a neighborhood of $M$ in $W$ is identified with a
neighborhood of $M$ in $M\times \R^k$ with the symplectic form
\eqref{eq:normal-form}.  Using this identification, we denote by $U_r$,
with $r>0$ sufficiently small, the neighborhood of $M$ in $W$
corresponding to $M\times D^k_r$, where $D^k_r$ is the ball of radius
$r$. Recall also that by definition $\rho=(p_1^2+\cdots+p_k^2)/2$, where
$(p_1,\ldots,p_k)$ are the coordinates on $\R^k$. (Thus
$U_r=\{\rho<r^2/2\}$.)

Let $K$ be a smooth function on $[0,r]$ such that
\begin{itemize}
\item $K$ is monotone decreasing and $K\equiv 0$ near $r$;
\item all odd-order derivatives of $K$ at $0$ are zero,
 and $K''(0)<0$ is close to zero.
\end{itemize}
Abusing notation, we also denote by $K$ the function on $W$ equal to
$K(|p|)$ on $U_r$, where $|p|=\sqrt{2\rho}$, and extended to be
identically zero outside $U_r$.  Then the Hamiltonian $K$ satisfies
the hypotheses of Proposition \ref{prop:orbits-deg}. By Proposition
\ref{prop:flow}, the Hamiltonian flow of $K$ on $U_r$ is a
reparametrization of the leaf-wise geodesic flow on $M$. Outside
$U_r$, the flow is the identity map.

\begin{Theorem}
\label{thm:main2}
Assume that a Hamiltonian diffeomorphism $\varphi_H$ displaces $M$
and that $r>0$ is sufficiently small so that, in particular,
$\varphi_H$ also displaces $U_r$. Then there exists a constant
$\Delta>0$, independent of $r$ and $K$, such that $K$ has a
contractible periodic orbit $\gamma$ with
$$
\Delta \leq A_K(\gamma)-\max K\leq \parallel H \parallel,
$$
provided that $\max K$ is large enough.
\end{Theorem}

\begin{proof}
The proof of the theorem relies on Proposition \ref{prop:orbits-deg} and
the following two lemmas, which are essentially contained in \cite{Bo1,Bo2}
and which will also be used in the proof of Theorem \ref{thm:main}.

For a closed curve $\eta\colon S^1\to M$ lying in a leaf of the foliation
$\FF$, set
$$
\delta(\eta)=\sum_{i=1}^k\left|\int_\eta\alpha_i\right|.
$$

\begin{Lemma}
\label{lemma:lower-bound}
There exists a constant $\delta_M>0$ such that
$$
\delta_M\leq \delta(\eta)
$$
for all non-trivial closed geodesics $\eta$ of the leaf-wise metric
$\alpha_1^2+\cdots+\alpha^2_k$ (see Proposition \ref{prop:flow}).
\end{Lemma}

\begin{proof}
Note that since, by Proposition \ref{prop:flow}, the metric is leaf-wise flat
we have $\alpha_i(\dot{\eta}(t))=\const$ for every leaf-wise geodesic $\eta$.
Hence
$$
\left|\int_\eta\alpha_i\right|=\int_{S^1} |\alpha_i(\dot{\eta}(t))|\,dt.
$$
As a consequence,
\begin{equation}
\label{eq:length}
\delta(\eta)\geq\length(\eta)
\end{equation}
which follows immediately from the fact that
$\sum_{i}|\alpha_i(v)|\geq 1$ for every unit vector $v$ tangent to $\FF$.

Assume that the required constant $\delta_M>0$ does not exist, i.e.,
$\delta(\eta_j)\to 0$ for some sequence of closed non-trivial
leaf-wise geodesics $\eta_j$. Thus, we also have $\length(\eta_j)\to
0$. Passing if necessary to a subsequence, we conclude that all
$\eta_j$ are contained in a small neighborhood in $M$.  (Indeed
extending the metric from $\FF$ to $M$, we can view $\eta_j$ as a
sequence of closed curves whose length converges to zero. Now it is
clear that the sequence $\eta_j$ contains a subsequence lying in a
neighborhood of arbitrarily small radius.)  Without loss of generality
we may assume that this neighborhood is foliated. Then, when $j$ is
large enough, $\eta_j$ is contained in a small ball in the leaf
$F\supset \eta_j(S^1)$ of $\FF$.  In particular $\eta_j$ is
contractible in $F$.  This is impossible since the metric is leaf-wise
flat.
\end{proof}

\begin{Remark}
When $M$ has restricted contact type, 
Lemma \ref{lemma:lower-bound} implies that the ``area
spectrum'' of $M$ is separated from zero by $\delta_M/k$, for
$\delta(\eta)=k|A(\eta)|$.
\end{Remark}

Fix a neighborhood $V=U_R$, for some $R>0$. In what follows, we will
always assume that $0<r<R/2$ and that $U_r$ is displaced by $\varphi_H$.

\begin{Lemma}
\label{lemma:low-bound2}
There exists a constant $c_V>0$, depending on $V$ but not on $K$ or
$r>0$, such that for any Floer anti-gradient trajectory $u$ for $K$
partially asymptotic to a non-trivial one-periodic orbit $\gamma$ at
$-\infty$ and to a point of $M$ at $\infty$, we have
$$
c_V \delta (\pi(\gamma))\leq A_K(\gamma)-\max K.
$$
\end{Lemma}

\begin{Remark}
Note that the right hand side of this inequality is the energy $E(u)$.
\end{Remark}

\begin{proof}
Let $f$ be a non-negative, smooth, decreasing function on $[0,R]$
identically equal to $1$ on $[0,R/2]$ and vanishing near $R$. Abusing
notation, we also denote by $f$ the function on $W$ equal to $f(|p|)$
on $V=U_R$, where $|p|=\sqrt{2\rho}$ as above, and extended to be
identically zero outside $V$.

Following \cite{Bo2}, set $\beta_i=f\pi^*\alpha_i$. This is a smooth
one-form on $W$. A feature of the form $\beta_i$, important in what
follows, is that

\begin{equation}
\label{eq:0}
i_{X_K}d\beta_i=0.
\end{equation}
To prove this, we first note that \eqref{eq:0} trivially holds outside
$U_r$. On the other hand, $f\equiv 1$ on $U_r$ since $r<R/2$. Thus
on $U_r$ we have $\beta_i=\pi^*\alpha_i$ and, point-wise,
$$
i_{X_K}d\beta_i=i_{X_K}\pi^*d\alpha_i=K'(\rho)i_{\pi_*X_\rho}d\alpha_i=0.
$$
The last equality follows from the fact that $\pi_*X_\rho$ is tangent to
$\FF$ by Proposition~\ref{prop:flow} and that $T\FF \subset \ker d\alpha_i$
since $M$ is stable.

Set $c_i=\parallel d\beta_i\parallel_{C^0}>0$ so that
$$
|d\beta_i(X,Y)|\leq c_i\parallel X\parallel\cdot \parallel Y\parallel,
$$
for any two tangent vectors $X$ and $Y$. Here, on the right hand side,
the norm is taken with respect an arbitrary metric compatible with
$\omega$ as in Section \ref{sec:statement}. (We emphasize that, since
the argument relies only on the results of Section \ref{sec:deg}, which
hold for general Hamiltonians and metrics, the metric need not meet
any regularity requirements for Floer anti-gradient trajectories; cf.
Remark \ref{rmk:orbits-deg}.) It is clear that $c_i$ is independent of
$K$.

We claim that
\begin{equation}
\labell{eq:energy-bound}
A_K(\gamma)-\max K=E(u)
\geq c_i^{-1}\left|\int_{\pi(\gamma)}\alpha_i\right|.
\end{equation}
The assertion of the lemma immediately follows from \eqref{eq:energy-bound}
by adding up these inequalities for $i=1,\ldots,k$ and setting
$c_V=\min c_i^{-1}/k$.

To prove \eqref{eq:energy-bound}, fix $s_j^\pm\to \pm\infty$ such that
$u(s_j^+)$ converges to a point of $M$ and $u(s_j^-)$ converges to $\gamma$
in $C^\infty(S^1,W)$. Then utilizing the definition of $c_i$ and \eqref{eq:0}
in the last step, we have
\begin{eqnarray*}
E(u)
&=&
\int_{\R\times S^1}
\left\| \frac{\p u}{\p s}\right\|
\cdot
\left\| \frac{\p u}{\p t}-X_K\right\|
\,ds\,dt
\\
&\geq&
c_i^{-1}
\int_{\R\times S^1}
\left|d\beta_i\left(\frac{\p u}{\p s}, \frac{\p u}{\p t}-X_K\right)\right|
\,ds\,dt
\\
&\geq&
c_i^{-1}
\lim_{j\to\infty}
\left|\int_{s_j^-}^{s_j^+}\int_{S^1}
d\beta_i\left(\frac{\p u}{\p s}, \frac{\p u}{\p t}-X_K\right)
\,dt\,ds\right|
\\
&=&
c_i^{-1}
\lim_{j\to\infty}
\left|\int_{s_j^-}^{s_j^+}\int_{S^1}
d\beta_i\left(\frac{\p u}{\p s}, \frac{\p u}{\p t}\right)
\,dt\,ds\right|.
\end{eqnarray*}
Note that at this stage we still do not know if the limit in question
exists. However, by applying Stokes' formula, we
see that
$$
\left|\int_{s_j^-}^{s_j^+}\int_{S^1}
d\beta_i\left(\frac{\p u}{\p s}, \frac{\p u}{\p t}\right)
\,dt\,ds\right|
=
\left|\int_{u(s_j^+)}\beta_i-\int_{u(s_j^-)}\beta_i\right|
\longrightarrow
\left|\int_{\gamma}\beta_i\right|
$$
as $j\to\infty$.
Furthermore, recall that $\gamma$ is contained in $U_r$
and $f|_{U_r}\equiv 1$. Thus,
$$
\left|\int_{\gamma}\beta_i\right|=\left|\int_{\gamma}\pi^*\alpha_i\right|
=\left|\int_{\pi(\gamma)}\alpha_i\right|,
$$
which completes the proof of \eqref{eq:energy-bound} and of the lemma.
\end{proof}

To finish the proof of the theorem, we apply Proposition \ref{prop:orbits-deg}
to the Hamiltonian $K$. Next, applying Lemma \ref{lemma:low-bound2} to
$u$ and $\gamma$ whose existence is guaranteed by this proposition, we have
\begin{equation}
\labell{eq:main-ineq}
\| H\|\geq A_K(\gamma)-\max K\geq c_V\delta(\pi(\gamma))\geq c_V \delta_M =:
\Delta>0,
\end{equation}
where the last inequality follows from Lemma \ref{lemma:lower-bound}
and the fact that $\pi(\gamma)$ is, up to parametrization and orientation,
a closed non-trivial geodesic of the leaf-wise geodesic flow on $\FF$.
\end{proof}

\begin{proof}[Proof of Theorem \ref{thm:main}]
Let $K$ be as above and let $\max K$ be large enough. Assertion~(i)
follows immediately from Theorem \ref{thm:main2}.

To prove (ii), we argue as in \cite{Bo2}.  Pick
$\gamma$ as in the proof of Theorem \ref{thm:main2}.  This is a
non-trivial periodic orbit and, since $K$ attains its
maximum along $M$, we see from \eqref{eq:main-ineq} that
\begin{equation}
\label{eq:A}
A(\gamma)\geq A_K(\gamma)-\max K \geq c_V\delta(\pi(\gamma)).
\end{equation}
Then, utilizing the normal form \eqref{eq:normal-form} and setting
$\eta=\pi(\gamma)$, we have
\begin{eqnarray*}
A(\eta)&=& A(\gamma)-\int_\gamma\sum p_i\alpha_i\\
&\geq& A(\gamma)-r\sum\left|\int_{\eta}\alpha_i\right|\\
&\geq& c_V\delta(\eta)-r\delta(\eta)\\
&=& (c_V-r)\delta(\eta)\\
&\geq& (c_V-r)\delta_M>0,
\end{eqnarray*}
provided that $r>0$ is small enough. This proves (ii)
for $\eta=\pi(\gamma)$. Moreover, we have also
established the \emph{a priori} lower bound
\begin{equation}
\label{eq:lb}
A(\eta)\geq (c_V-r)\delta_M .
\end{equation}

It remains to prove (iii). To this end, we first need to establish a general
property of the metric $\rho$, which holds for arbitrary stable
coisotropic submanifolds.  Recall that the length spectrum of a
(leaf-wise) metric on a foliated manifold $(M,\FF)$ is the collection
of lengths of all closed leaf-wise geodesics in $M$. In contrast with
the length spectrum of a metric on $M$, the length spectrum of a
leaf-wise metric need not in general be nowhere dense, even if the
metric is leaf-wise flat.  However, as our next observation shows,
this is the case for the metric $\rho$.

\begin{Lemma}
\label{lemma:nowheredense}
Assume that $M$ is stable. The length spectrum of $\rho$ is closed, has
zero measure, and is, therefore, nowhere dense.
\end{Lemma}

\begin{proof}
The assertion that the length spectrum is closed holds obviously for
any metric. To show that the spectrum has zero measure, consider a
level $\{\rho=\const\}\subset M\times\R^k$ and recall that by
Proposition \ref{prop:flow} closed geodesics of $\rho$ are
projections to $M$ of closed characteristics on this level.
Furthermore, again by Proposition \ref{prop:flow}, the length of a
geodesic is equal to the integral of $\lambda=\sum p_i\alpha_i$ over
the corresponding closed characteristics. (Strictly speaking, the
equality holds only up to a factor depending only on the level, and
we chose $\const$ so that this factor is equal to one.) The
characteristic foliation on the level is tangent to the distribution
$\ker d\lambda$ as is easy to see from \eqref{eq:spray}.
It follows that the integrals of $\lambda$ over
closed characteristics form a zero measure set.  (The distribution
$\ker d\lambda$ need not have constant rank, but the standard
argument still applies.)
\end{proof}

Let now $\gamma$ be as above and $\eta=\pi(\gamma)$. Our next goal is
to obtain an upper bound on $A(\gamma)$ and $A(\eta)$. To this end, we
need to impose an extra requirement on $K$ guaranteeing that
$K(\gamma)$ is close to $\max K$.

Namely, fix $r_-$ and $r_+$ such that $0<r_-<r_+<r$ and pick a sufficiently
small constant $\eps>0$.  We can chose $K$
so that on $U_r$, we have
\begin{itemize}

\item $\max K-K<\eps$ on $[0,r_-]$ and $K<\eps$ on $[r_+,r]$;

\item on $[r_-,r_+]$ the Hamiltonian $K$,  thought of
as function of $|p|$, has constant slope
lying outside the length spectrum of the metric $\rho$;

\item $K$ has a sufficiently large variation over the constant slope
range, i.e., $K(r_-)-K(r_+)> C\cdot\| H\|$, where $C$ is a constant,
to be specified later, depending only on $V$.

\end{itemize}
In particular, these conditions ensure that the only non-trivial periodic
orbits of $K$ occur within the shells $(0,r_-)$ and $(r_+,r)$.

Furthermore, we claim that $\gamma$ lies in the shell $(0,r_-)$, and hence
\begin{equation}
\label{eq:nearmax}
\max K-K(\gamma)<\eps.
\end{equation}
Indeed, first observe that since $M$ has restricted contact type,
we have
$$
A(\eta)=\int_\eta\alpha_i\text{ for } i=1,\ldots,k.
$$
Thus, by \eqref{eq:energy-bound},
$$
|A(\eta)|\leq c_i\| H\|\text{ for } i=1,\ldots,k
$$
and, as a consequence,
\begin{eqnarray*}
|A(\gamma)| &=& \left|A(\eta)+\int_\gamma \sum_i p_i\alpha_i\right|\\
           &\leq& |A(\eta)|+r\sum_i\left|\int_\eta\alpha_i\right|    \\
           &=& (1+kr)|A(\eta)|\\
           &\leq& C\cdot\|H\|,
\end{eqnarray*}
where, for instance, $C=(1+kr)c_1$.
On the other hand, if $\gamma$ were in the shell $(r_+,r)$, we
would have
$$
\max K- K(\gamma)> K(r_-)-K(r_+)> C\cdot\| H\|,
$$
and hence,
$$
A(\gamma)\geq\max K- K(\gamma)> C\cdot\| H\|.
$$
Thus $\gamma$ is indeed in the shell $(0,r_-)$ and \eqref{eq:nearmax} holds.

Next note that
$$
A(\eta)=\big(A_K(\gamma)-\max K\big)
+\big(\max K-K(\gamma)\big) -\int_\gamma\sum p_i\alpha_i.
$$
Here the first term in bounded from above by $\| H \|$ due to Proposition
\ref{prop:orbits-deg}, the second term is bounded from above by $\eps$ due to
\eqref{eq:nearmax}, and the last term is bounded from above by
$r\delta(\eta)\leq r c_V^{-1} \| H\|$ according to \eqref{eq:main-ineq}.
As a consequence,
\begin{equation}
\label{eq:hplast}
(c_V-r)\delta_M\leq A(\eta)\leq\| H\| +\eps + rc_V^{-1}\| H\|,
\end{equation}
where the first inequality is \eqref{eq:lb}.
Note that we also have
an \emph{a priori} upper bound on the length of $\eta$,
$$
\length(\eta)\leq \delta(\eta)\leq c_V^{-1}\| H\|,
$$
resulting from \eqref{eq:length} and \eqref{eq:main-ineq}.

Consider a sequence of Hamiltonians $K_i$ as above with $r_i\to 0$ and
$\eps_i\to 0$. Then by the Arzela--Ascoli theorem the geodesics $\eta_i$
converge to a geodesic $\eta$ and, passing to the limit in
\eqref{eq:hplast}, we see that
$$
0<c_V\delta_M\leq A(\eta)\leq\| H\|,
$$
which concludes the proof.
\end{proof}

\begin{Remark}
\labell{rmk:Dr} The first assertion of Theorem \ref{thm:main} can also
be established by making use of Proposition \ref{prop:orbits-deg2}
instead of Proposition \ref{prop:orbits-deg}, although in a somewhat
less direct way. Let $K$ and $r$ be as in Theorem \ref{thm:main2} with
$\max K>\| H\|$ and let $\eps>0$ be so small that $f=\eps K$ is
$C^2$-small. Denote by $u$ and $\gamma$ the homotopy connecting
trajectory and the periodic orbit from Proposition
\ref{prop:orbits-deg2} for a \emph{linear} homotopy from $K$ to $f$.
Since $\max K>\| H\|> \sigma(K)$, we conclude that the orbit $\gamma$
is non-trivial. It is not hard to see that the proof of Lemma
\ref{lemma:low-bound2} goes through for $u$. (Note that the
assumptions that $f=\eps K$ and that the homotopy is linear are
essential to make sure that \eqref{eq:0} holds for every Hamiltonian
$K^s$ in the homotopy.)  Hence, as in the proof of Theorem
\ref{thm:main2}, we have
$$
\|H\|\geq E(u)\geq c_V\delta(\pi(\gamma))\geq c_V\delta_M,
$$
which proves assertion (i). Even though this approach
does not lead to assertion (ii) in its full generality, it does imply that
$\FF$ has a leaf-wise non-trivial closed geodesic, contractible in $W$.
This is an analogue of the existence result for closed characteristics
on closed, stable or contact type hypersurfaces, \cite{HZ,vi}; see also
the survey \cite{gi:alan} for a discussion of more recent results.

Finally note that passing to the limit as $f\to 0+$, we infer from
Proposition \ref{prop:orbits-deg2} that $\sigma(K)\geq c_V\delta_M$.
As a consequence, $\chom(U)\geq c_V\delta_M$ for any neighborhood
$U$ of $M$. Thus $\chom(M)\geq c_V\delta_M>0$, whenever $M$ is stable
and displaceable; cf.\ \cite{Dr}.
\end{Remark}

\begin{proof}[Proof of Theorem \ref{thm:comm}]
It suffices to show that the levels $M_a$ carrying an orbit $\gamma$
with the required properties exist arbitrarily close to $M$.  Set
$K=f(K_1,\ldots,K_k)$, where $f\colon\R^k\to\R$ is a bump-function
supported in a small neighborhood of the origin in $\R^k$ and such
that $\max f$ is large enough. Since the support of $f$ is small, we
may assume that $\supp K$ is displaceable and all $a\in \supp f$ are
regular values of $\vec{K}$.  By Proposition \ref{prop:orbits-deg},
the flow of $K$ has a contractible in $W$ one-periodic orbit
$\gamma$ with $A(\gamma)>0$; see, e.g., Remark \ref{rmk:orbits-deg}.
Furthermore, the Hamiltonian $K$ Poisson--commutes with all $K_i$, and
thus $\gamma$ is tangent to a regular level $M_a$.
\end{proof}

\section{Leaf-wise intersections for hypersurfaces}
\label{sec:leafwise}

The goal of this section is to prove the leaf-wise intersection
property (Theorem~\ref{thm:leafwise}) for hypersurfaces of
restricted contact type in subcritical Stein manifolds.
In fact the theorem holds for a slightly broader class of
ambient manifolds than subcritical Stein. Namely, let $W$ be an
exact symplectically aspherical manifold convex in the sense of
\cite{FS} and let as above $M$ be a closed hypersurface of restricted
contact type in $W$ bounding a domain $U$. Assume that $H$ is a compactly
supported Hamiltonian on $W$ such that
$$
\|H\|<\chom(U)
$$
and
\begin{equation}
\label{eq:upper-chom}
\chom\big(\bar{U}\cup\supp(H)\big)<\infty.
\end{equation}
Then Theorem \ref{thm:leafwise} (in a slightly generalized form)
asserts that $\varphi(F)\cap F\neq\emptyset$ for some leaf $F$ of the
characteristic foliation of $M$ and $\varphi=\varphi_H$.

It is worth mentioning that the requirement \eqref{eq:upper-chom} is
quite restrictive.  (Note that \eqref{eq:upper-chom} implies that
$\chom(U)<\infty$.)  The reason that \eqref{eq:upper-chom} holds
if $W$ is a subcritical Stein manifold is that every compact set in $W$
is displaceable; see, e.g., \cite{BC}. However, \eqref{eq:upper-chom}
fails, as can be seen from the results of \cite{vi:functors}, when $U$
is a tubular neighborhood of the zero section in a cotangent bundle.
Essentially the only case where \eqref{eq:upper-chom} can be verified
is where $\bar{U}\cup\supp(H)$ is displaceable in $W$.

Hereafter, we assume that $W$ and $M$ are as above. In this case, the
action selector $\sigma$ from Section \ref{sec:sel-open} can also be
defined as in \cite{FS} for all compactly supported Hamiltonians (not
necessarily positive) and has, in addition to (S0)--(S4), a number of
other properties. For instance, $\sigma(H)=0$, whenever $H\leq 0$, and
$\sigma$ is sub-additive, i.e., $\sigma(H\# K)\leq
\sigma(H)+\sigma(K)$; see \cite{FS}.  Furthermore, $\sigma(K)$ depends
only on the time-one flow $\varphi_K$ and we will also use the
notation $\sigma(\varphi_K)$.

\begin{proof}[Proof of Theorem \ref{thm:leafwise}]
First note that as in the proof of Proposition
\ref{prop:orbits-nondeg} we may require that $\min H_t=0$ for all $t$
and hence $\| H\|= E^+(H)$ and $H\geq 0$. (This can be achieved,
keeping \eqref{eq:upper-chom}, by replacing $H$ by the Hamiltonian
$f\cdot (H-\min H)$, where $f$ is a cut-off function identically equal
to one near $\bar{U}\cup\supp(H)$.)  Furthermore, without loss of
generality we may also assume that the fixed point set $\Fix(\varphi)$
does not meet $M$, for otherwise the assertion is obvious. Then
$\varphi$ has no fixed points near $M$.

Fix a global primitive $\lambda$ of $\omega$, restricting to a contact
form on $M$. Let $U_\eps=M\times [-\eps,0]$ be a narrow shell inside
$U$, containing no points of $\Fix(\varphi)$.
We choose the projection $U_\eps\to
[-\eps,0]$ so that the Hamiltonian flow of this projection (thought of
as a function on $U_\eps$) is exactly equal to the Reeb flow of
$\lambda|_{M_\tau}$, where $M_\tau=M\times\{\tau\}$, for all $\tau\in
[-\eps,0]$. The projection $U_\eps\to M$ is chosen so that the
characteristic foliation on $M_\tau$ projects to the characteristic
foliation on $M$. (It is well known that such a shell exists, when $M$
has contact type.)

Let $K\geq 0$ be a non-negative function which is equal to zero on
$W\ssminus U$ and to $\max K>0$ on $U\ssminus
U_\eps$, and which is a monotone decreasing function of $\tau$ on
$U_\eps$.
Set $\psi_s=\varphi\varphi_K^s=\varphi\varphi_{sK}$ for $s\in [0,\,1]$.
Then we have the following disjoint union decomposition
$$
\Fix(\psi_s)=\Fix(\varphi)\sqcup Z_s, \text{ where }
Z_s=\{x\in U_\eps \mid \varphi^{-1}(x)=\varphi_{sK}(x)\},
$$
as is easy to see recalling that $\Fix(\varphi)\cap U_\eps=\emptyset$.

Since $H\geq 0$, the Hamiltonian generating the map $\varphi^{-1}$ is
non-positive and hence $\sigma(\varphi^{-1})=0$. Then, by
conjugation invariance  and sub-additivity of the action selector
$\sigma$ (see \cite{FS} and also \cite{schwarz,Vi:gen}), we have
\begin{equation}
\label{eq:long-sigma}
\sigma(\varphi_{sK})=\sigma(\varphi\varphi_{sK}\varphi^{-1})
\leq \sigma(\varphi\varphi_{sK}) + \sigma(\varphi^{-1})=
\sigma(\varphi\varphi_{sK})=\sigma(\psi_s).
\end{equation}

From the monotonicity property, (S0), of $\sigma$ we infer
that $\chom(U)=\sup\sigma(\varphi_K)$, where the supremum is taken over
all $K$ as above.
Pick $K$ such that $\sigma(\varphi_K)$ is close to $\chom(U)$:
\begin{equation}
\label{eq:upperH}
E^+(H)<\sigma(\varphi_K)\leq \chom(U).
\end{equation}
As $s$ varies through the interval $[0,\,1]$, the action selector
$\sigma(\psi_s)\in \CS(\psi_s)$
changes from
$$
\sigma(\psi_0)=\sigma(\varphi)\leq E^+(H)
$$
to
$$
\sigma(\psi_1)\geq \sigma(\varphi_K)> E^+(H),
$$
where the first inequality follows form \eqref{eq:long-sigma} and the
second one from \eqref{eq:upperH}. By continuity of the action
selector, we see that $\CS(\psi_s)$ cannot be independent of $s$.
Furthermore, $\psi_0=\varphi$, and thus $\Fix(\psi_0)=\Fix(\varphi)$ and
$Z_0=\emptyset$. As a
consequence, the part $Z_s$ of $\Fix(\psi_s)$ must be non-empty for
some $s_0\in (0,\,1]$ and, moreover, $\sigma(\psi_{s_0})$ is the
action value of $\psi_{s_0}$ on $x\in Z_{s_0}$.

As in the proof of Proposition \ref{prop:orbits-nondeg}, we may
reparametrize the Hamiltonians $H$ and $K$ (making $K$ now
time-dependent), without altering the time-one maps, the action spectra,
the Hofer norms, and the action selectors so that $K_t \equiv 0$ when $t\in
[1/2,\, 1]$ and $H_t \equiv 0$ when $t\in [0,\,1/2]$. From now on, we
assume that $H$ and $K$ have this property. We will also denote
$\varphi^t_H$ by $\varphi^t$.

Consider the orbit $\gamma(t)$ through $x=\gamma(0)$ of the
time-dependent flow $\psi_{s_0}^t=\varphi^t\varphi^t_{s_0K}$. Let
$G=H\#(s_0 K)$ be the Hamiltonian generating this flow. Due to the
above reparametrizations of $H$ and $K$, the orbit $\gamma$ is
comprised of two parts: $\gamma_1(t)=\varphi^t_{s_0K}(x)$ ending at
$y=\varphi_{s_0K}(x)$ and $\gamma_2(t)=\varphi^t(y)$ ending at
$x$. Note that $x$ and $y$ lie on the same Reeb orbit (i.e., a leaf of
characteristic foliation) on some level $M_\tau$ and
$\varphi(y)=x$. Furthermore,
\begin{equation}
\label{eq:gamma}
\sigma(\psi_{s_0})=A_G(\gamma)
=-\int_{\gamma_1}\lambda+ \int_0^{1/2} s_0K_t(\gamma_1(t))\,dt + A_H(\gamma_2),
\end{equation}
where
$$
A_H(\gamma_2)=
-\int_{\gamma_2}\lambda+ \int_{1/2}^1 H_t(\gamma_2(t))\,dt.
$$
The term $T=-\int_{\gamma_1}\lambda$ is the time required for the Reeb flow
on $M_\tau$ to move $x$ to $y$ and our next goal is to establish an upper
bound on $T$ independent of $K$.

Consider the function $f(z)$ equal to the action of $H$ on the orbit
$\varphi^t(z)$, $t\in [0,\,1]$, defined using the primitive $\lambda$.
(For instance, $f(y)=A_H(\gamma_2)$.) This function is independent of
$K$.  Clearly $f$ is a compactly supported function and $C=-\min f$ is
also independent of $K$ and $A_H(\gamma_2)\geq -C$. Furthermore, the
middle term in \eqref{eq:gamma} is non-negative, for $K\geq 0$. As a
consequence,
$$
\sigma(\psi_{s_0})\geq T-C.
$$
Finally note that the Hamiltonian $G$ generating $\psi_{s_0}$ is supported
in $U\cup\supp(H)$. Thus $\sigma(\psi_{s_0})\leq\chom\big(U\cup\supp(H)\big)$.
Therefore,
$$
T\leq \chom\big(U\cup\supp(H)\big)+ C\leq e\big(U\cup\supp(H)\big)+ C
<\infty,
$$
where the upper bounds on the right hand side are clearly independent of $K$.

To finish the proof, consider a sequence of Hamiltonians $K_i$ such as
$K$, non-constant on a more and more narrow range of $\tau$ in
$[-\eps,\,0]$ eventually converging to zero. For each $K_i$
we have a pair of points $x_i$ and $y_i$ lying on the same
Reeb orbit on some $M_{\tau_i}$ with $\tau_i\to 0$ and such
that $\varphi(y_i)=x_i$. Furthermore, the Reeb flow requires time
$T_i\leq \chom\big(U\cup\supp(H)\big)+ C$ to move $x_i$ to $y_i$. Applying
the Arzela--Ascoli theorem and passing if necessary to a subsequence,
we obtain points $x=\lim x_i$ and $y=\lim y_i$ on $M$ lying on the same
Reeb orbit and such that $\varphi(y)=x$.
This completes the proof of the theorem.
\end{proof}

\begin{Remark}
The proof of Theorem \ref{thm:leafwise} also yields the upper bound
$\chom\big(\bar{U}\cup\supp(H)\big)+ C$ for the ``Reeb distance'' from
$x$ to $y$.
\end{Remark}

In the next example we show that Theorem \ref{thm:leafwise}, as stated
with the upper bound on $\|H\|$, does not extend to hypersurfaces in
$\R^{2n}$ that do not have contact type. To be more precise, we
construct a Hamiltonian flow $\varphi^t$ on $\R^{2n}$ and a sequence
of hypersurfaces $M_i$, $C^0$-converging to $S^{2n-1}$,
such that $M_i$ and $\varphi^{t_i}(M_i)$ have no leaf-wise
intersections for some sequence of times $t_i\to 0+$.

\begin{Example}
\label{ex:leafwise-fail}
Let $S^{2n-1}$ be the unit sphere in
$\R^{2n}$ and let $\varphi^t$ be the Hamiltonian flow of $H=f\cdot p_1$,
where $(p_1,q_1,\ldots,p_n,q_n)$ are the standard coordinates on
$\R^{2n}$ and $f$ is a cut-off function equal to one near $S^{2n-1}$.
For $t>0$ small, the only leaf-wise intersections of
$S^{2n-1}$ and $\varphi^t(S^{2n-1})$ are two points $x_1$ and $x_2$
on the unit circle $S$ in the $(p_1,q_1)$-plane. (The points
$\varphi^t(x_1)$ and $\varphi^t(x_2)$ are the intersections of $S$
and the transported circle $S+(t,0)$.) Let us now insert two symplectic
plugs into $S^{2n-1}$ centered at points on $S$ between $x_1$ and
$\varphi^t(x_1)$ and between $x_2$ and $\varphi^t(x_2)$; see, e.g.,
\cite{Ci,Gi95,gi:bayarea,gg:ex2,ke:example}. We choose the plugs so
small and center them in such a way that they are displaced by
$\varphi^t$. As a result, we obtain a new hypersurface $M$ that is
$C^0$-close to $S^{2n-1}$, differs from $S^{2n-1}$ only within the
plugs, and such that the leaf $S$ is broken into two leafs: one
containing $x_1$ and $x_2$ and the other one containing
$\varphi^t(x_1)$ and $\varphi^t(x_2)$.  We claim that $M$ and
$\varphi^t(M)$ have no leaf-wise intersections.  Indeed, $x_1$ and
$x_2$ are no longer leaf-wise intersections for $M$ and
$\varphi^t(M)$, and since the plugs are displaced and due to the
plug-symmetry conditions, no new leaf-wise intersections are created.

Applying this construction to a sequence $t_i\to 0+$, we obtain a sequence
of perturbations $M_i$ of $S^{2n-1}$, $C^0$-converging to $S^{2n-1}$, and such
that $M_i$ and $\varphi^{t_i}(M_i)$ have no leaf-wise intersections. Note that
$\varphi^{t_i}\stackrel{C^\infty}{\to}\id$ while $\chom(U_i)\to \chom(U)=\pi$,
where $U_i$ is the domain bounded by $M_i$ and $U$ is the unit ball. It is
also clear that $\|H_i\|\to 0$, where $H_i=t_i H$ is a Hamiltonian generating
$\varphi^{t_i}$.
\end{Example}

\end{document}